\documentclass[12pt,reqno]{article}

\usepackage[usenames]{color}
\usepackage{amssymb}
\usepackage{graphicx}
\usepackage{amscd}

\usepackage[colorlinks=true,
linkcolor=webgreen,
filecolor=webbrown,
citecolor=webgreen]{hyperref}

\definecolor{webgreen}{rgb}{0,.5,0}
\definecolor{webbrown}{rgb}{.6,0,0}

\usepackage{color}
\usepackage{fullpage}
\usepackage{float}

\usepackage{graphics,amsmath,amssymb}
\usepackage{amsthm}
\usepackage{amsfonts}
\usepackage{latexsym}
\usepackage{epsf}

\usepackage{verbatim}

\usepackage{subfigure}
\usepackage{booktabs}
\usepackage{array}

\usepackage{tikz}
\usepackage{forest}

\usepackage{mathtools}

\theoremstyle{plain}

\theoremstyle{definition}

\theoremstyle{remark}

\usepackage{booktabs}
\usepackage{array}
\newcolumntype{M}[1]{>{\centering\arraybackslash}m{#1}}

\newcommand{\setsquarewithlabels}[1]{
        \def\type{#1}
        \setsquarecontinued
}
\newcommand{\setsquarecontinued}[9]{
\begin{center}
        \begin{tabular}{ccc}
            \includegraphics[width = 1.25cm]{#1.png} & \includegraphics[width = 1.25cm]{#2.png} & \includegraphics[width = 1.25cm]{#3.png}\\
            \includegraphics[width = 1.25cm]{#4.png} & \includegraphics[width = 1.25cm]{#5.png} & \includegraphics[width = 1.25cm]{#6.png}\\
            \includegraphics[width = 1.25cm]{#7.png} & \includegraphics[width = 1.25cm]{#8.png} & \includegraphics[width = 1.25cm]{#9.png}\\
        \end{tabular}
        \caption{Type \type.}
        \label{fig:\type}
\end{center}
}

\newcommand{\setsquare}[9]{
        \begin{tabular}{ccc}
            \includegraphics[width = 1.25cm]{#1.png} & \includegraphics[width = 1.25cm]{#2.png} & \includegraphics[width = 1.25cm]{#3.png}\\
            \includegraphics[width = 1.25cm]{#4.png} & \includegraphics[width = 1.25cm]{#5.png} & \includegraphics[width = 1.25cm]{#6.png}\\
            \includegraphics[width = 1.25cm]{#7.png} & \includegraphics[width = 1.25cm]{#8.png} & \includegraphics[width = 1.25cm]{#9.png}\\
        \end{tabular}
}

\begin{document}

\title{The Classification of Magic SET Squares}
\author{
Eric Chen \quad William Du \quad Tanmay Gupta \quad Alicia Li \\ Srikar Mallajosyula \quad Rohith Raghavan \quad Arkajyoti Sinha \\ Maya Smith \quad Matthew Qian \quad Samuel Wang\\
PRIMES STEP\\
\\
Tanya Khovanova\\
MIT
}

\maketitle

\begin{abstract}
A magic SET square is a 3 by 3 table of SET cards such that each row, column, diagonal, and anti-diagonal is a set. We allow the following transformations of the square: shuffling features, shuffling values within the features, rotations and reflections of the square.  Under these transformations, there are 21 types of magic SET squares. We calculate the number of squares of each type. In addition, we discuss a game of SET tic-tac-toe.
\end{abstract}

\section{The game of SET}

The game of SET is one of the most mathematical games. To play it you need a deck of special cards. Each card has one, two, or three identical objects drawn on it. The objects could be of three different colors: red, green, or purple. The object can have three different shapes: oval, diamond, and squiggly. The cards can also differ in shading. There are three types of shading: empty, full, and striped.

Overall, there are four features: number, shape, color, and shading. And there are three possibilities for each feature. The total number of cards in the deck is 81. Each card has a unique choice for one of the four features.

\begin{figure}[ht!]
        \begin{center}
        \begin{tabular}{ccc}
            \includegraphics[width = 1.25cm]{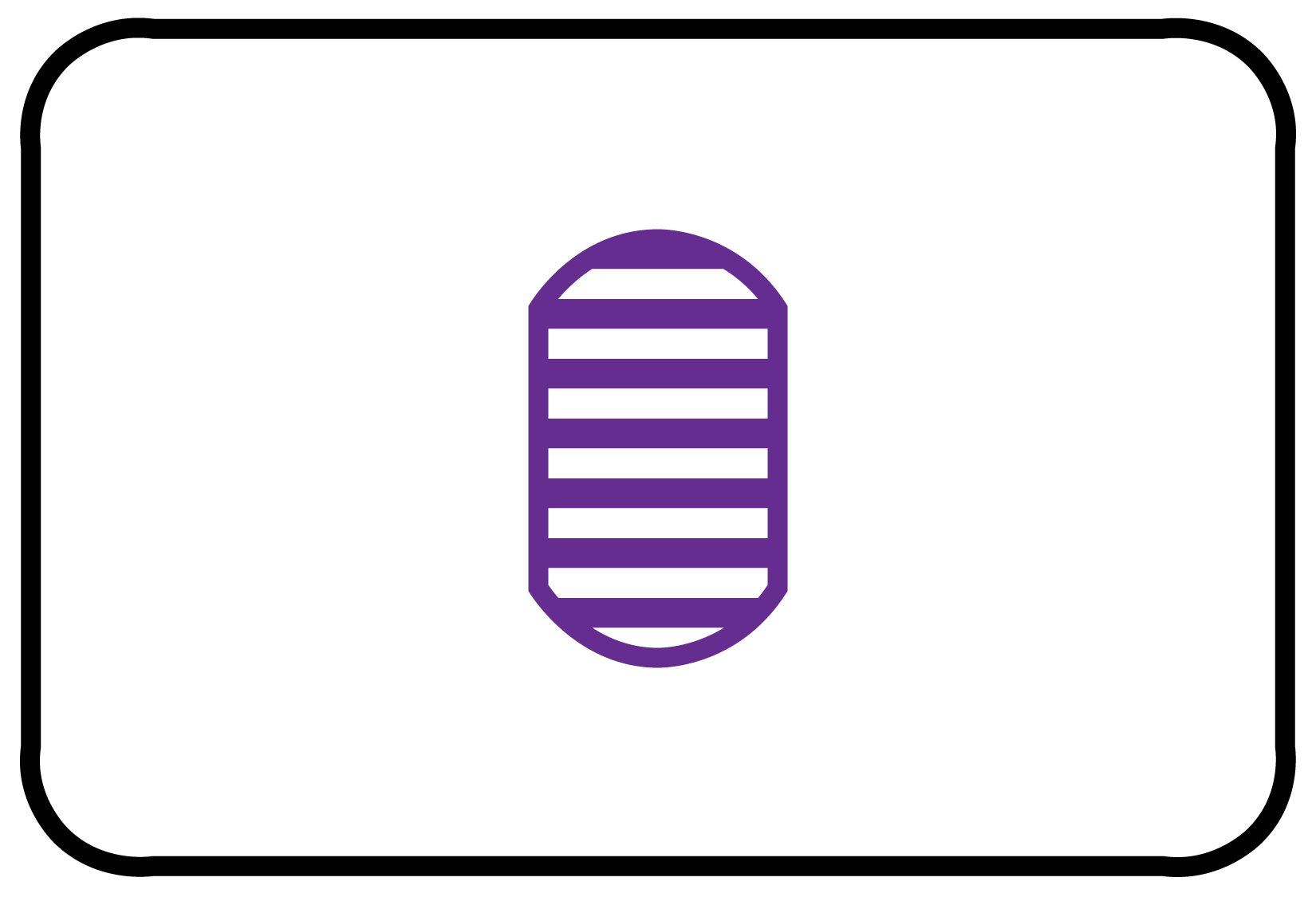} & \includegraphics[width = 1.25cm]{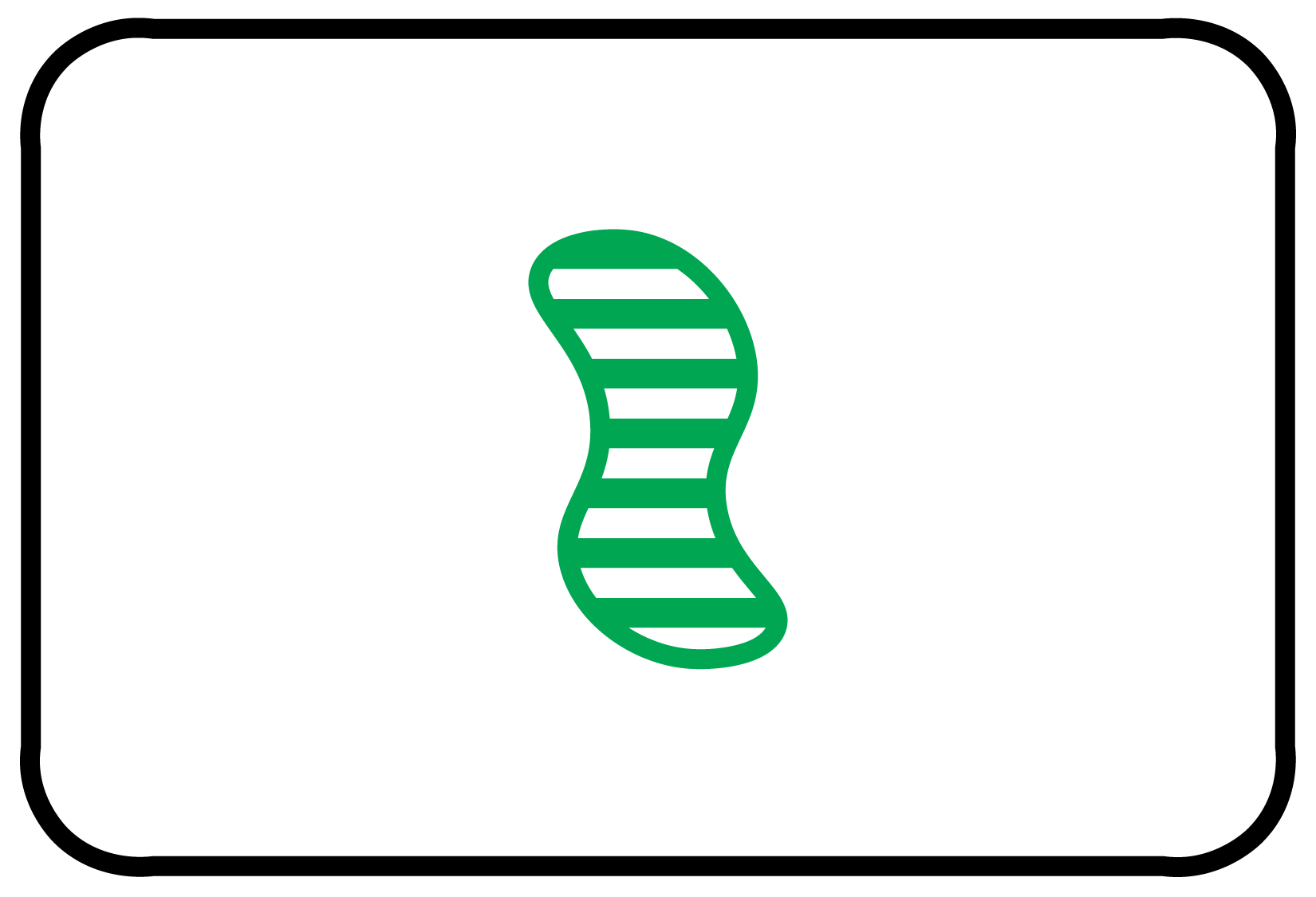} & \includegraphics[width = 1.25cm]{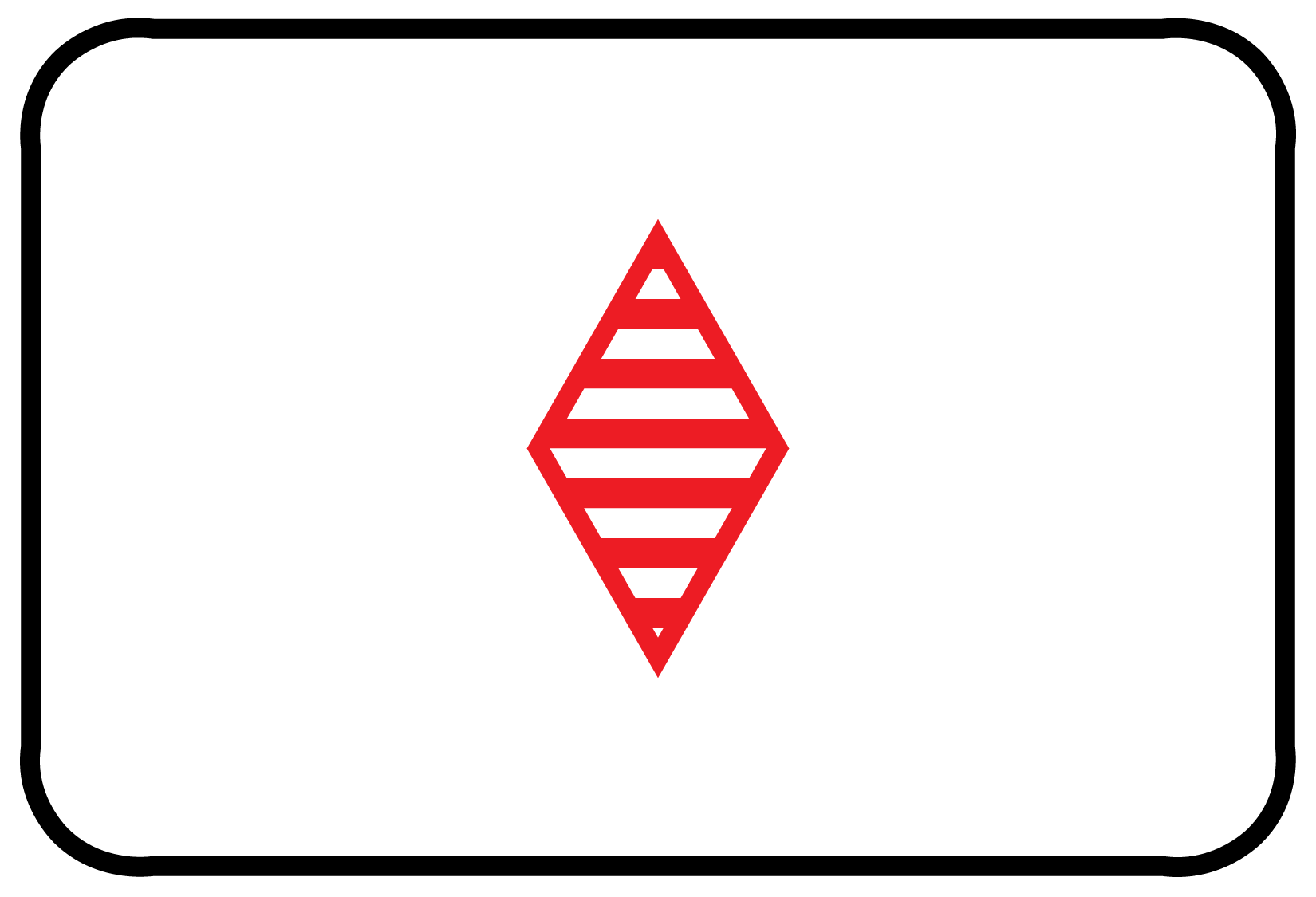}\\
        \end{tabular}
\caption{An example of a set.}
\label{fig:randomset}
        \end{center}
\end{figure}

Three cards can form \textit{a set} if and only if for every feature all three cards are either the same or all different. Figure~\ref{fig:randomset} shows an example of a set.

All three cards have the same number of objects---1, all the cards have the same shading---striped. All the cards have different colors and different shapes.

We say that set $A$ is of \textit{order} $k$ if there are exactly $k$ features that are different. The above example is a set of order 2. For many people, sets of order 1 are easier to see, and sets of order 4 are way more difficult to see.

\begin{figure}[!ht]
   \begin{minipage}{0.48\textwidth}
        \begin{tabular}{ccc}
            \includegraphics[width = 1.25cm]{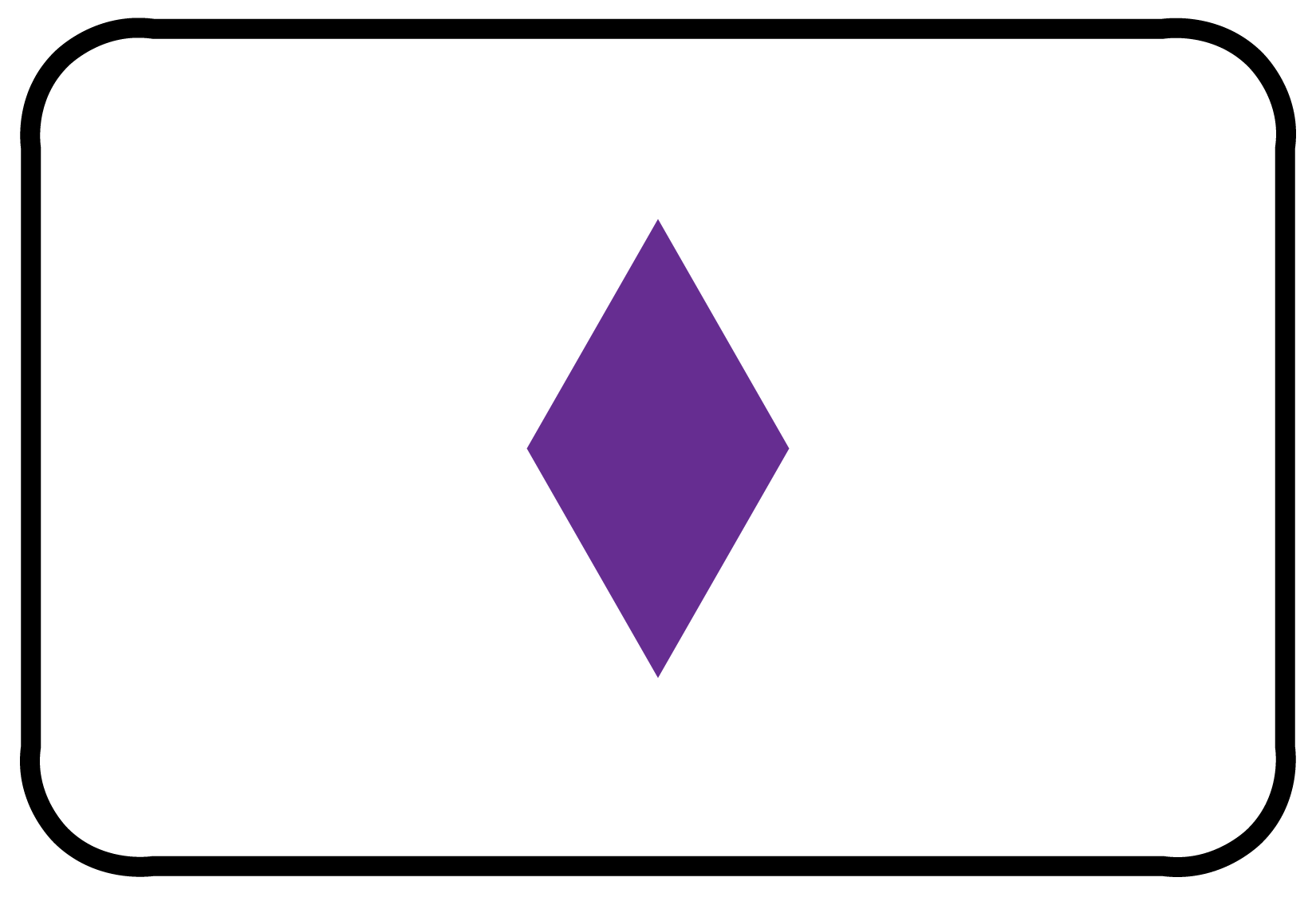} & \includegraphics[width = 1.25cm]{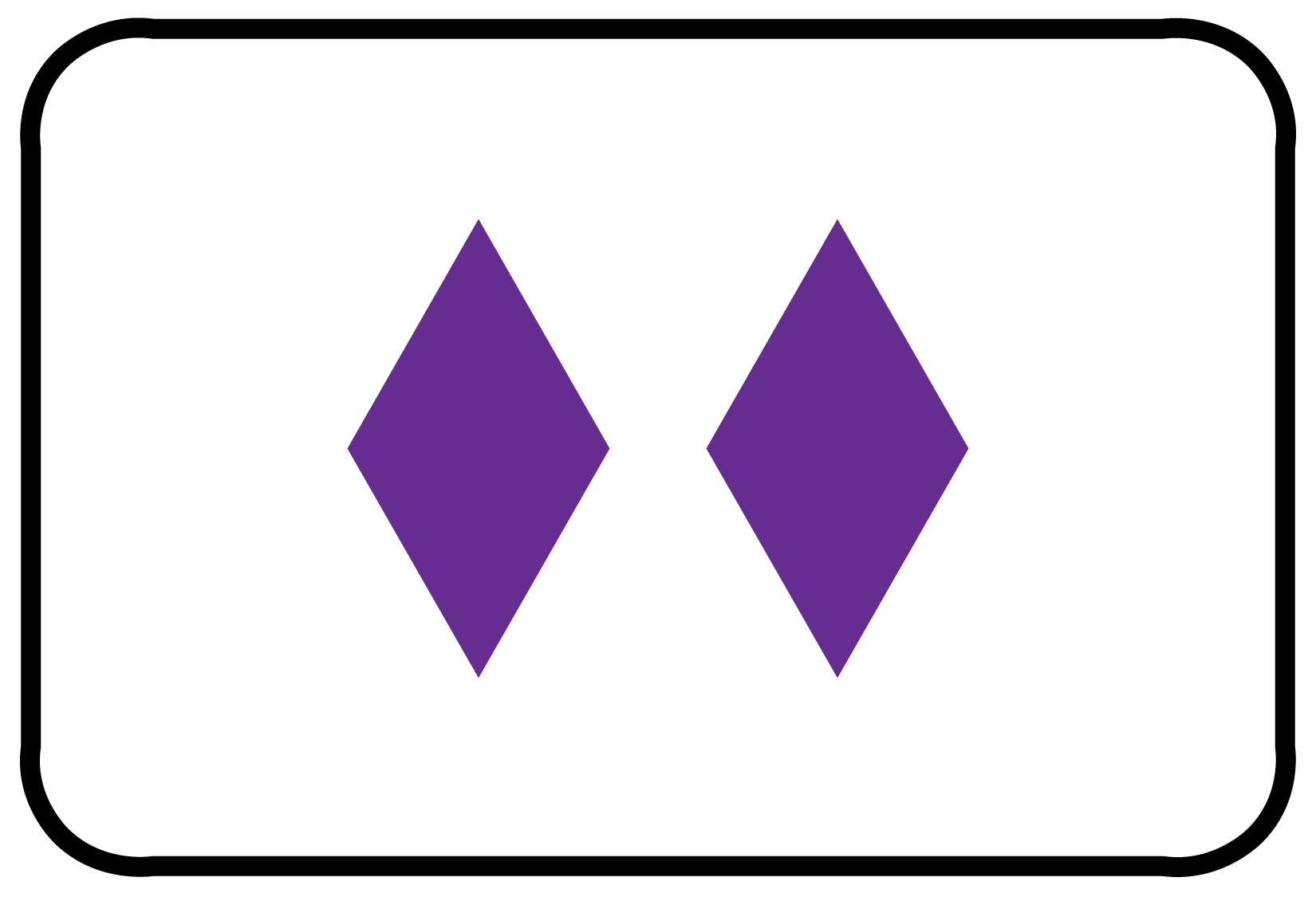} & \includegraphics[width = 1.25cm]{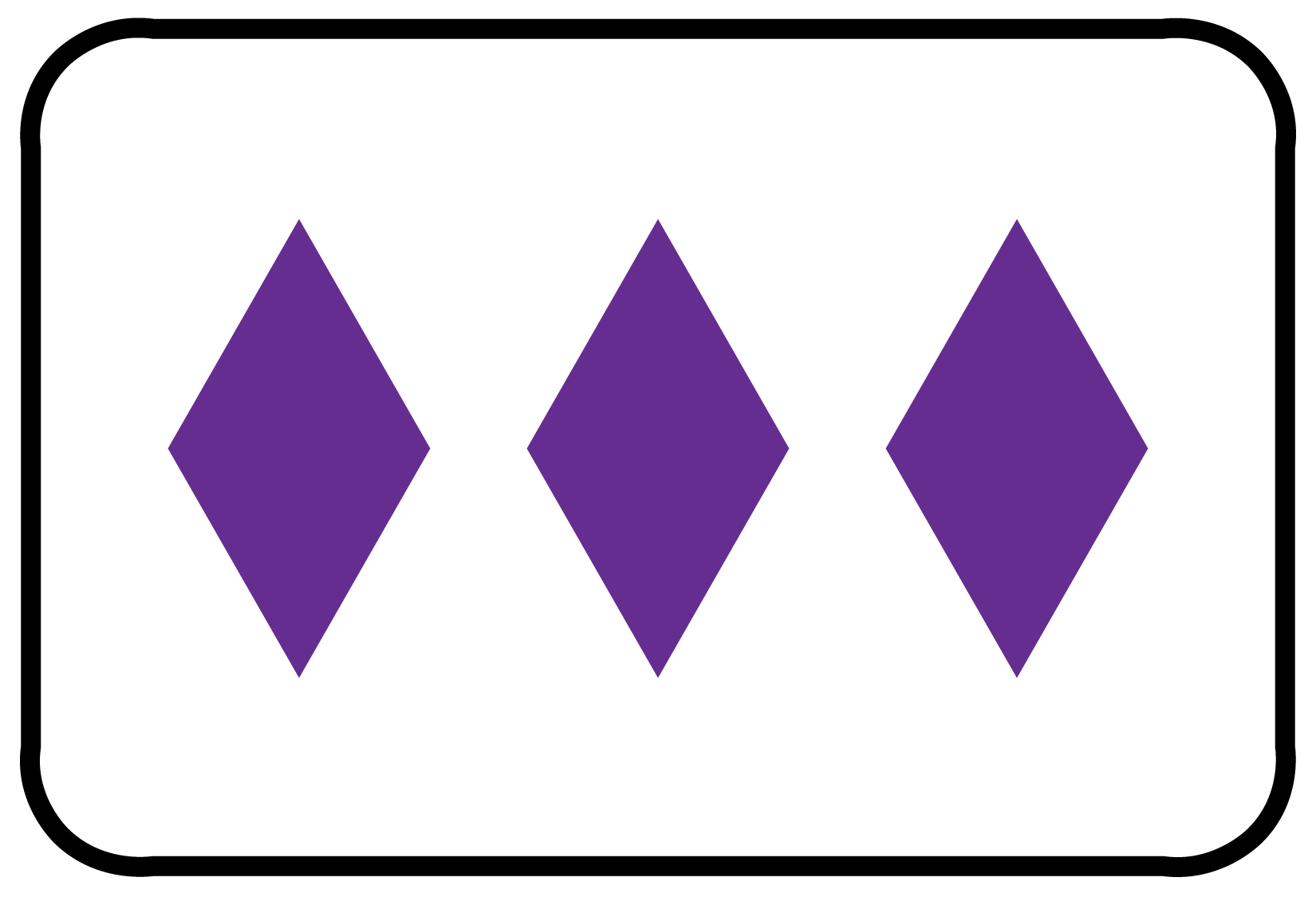}\\
        \end{tabular}
\caption{A set of order 1.}
\label{fig:set1}
   \end{minipage}\hfill
   \begin{minipage}{0.48\textwidth}
        \begin{tabular}{ccc}
            \includegraphics[width = 1.25cm]{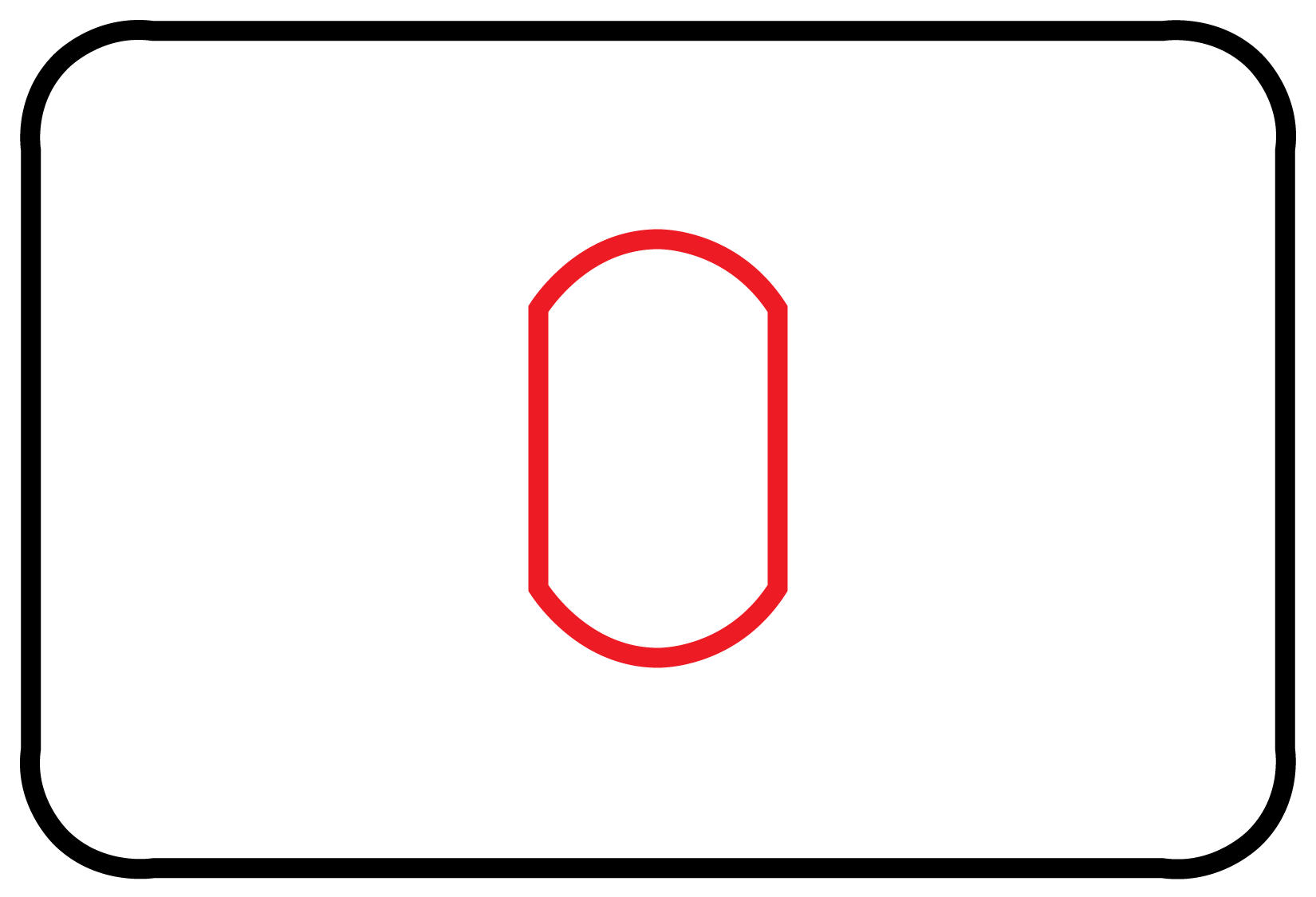} & \includegraphics[width = 1.25cm]{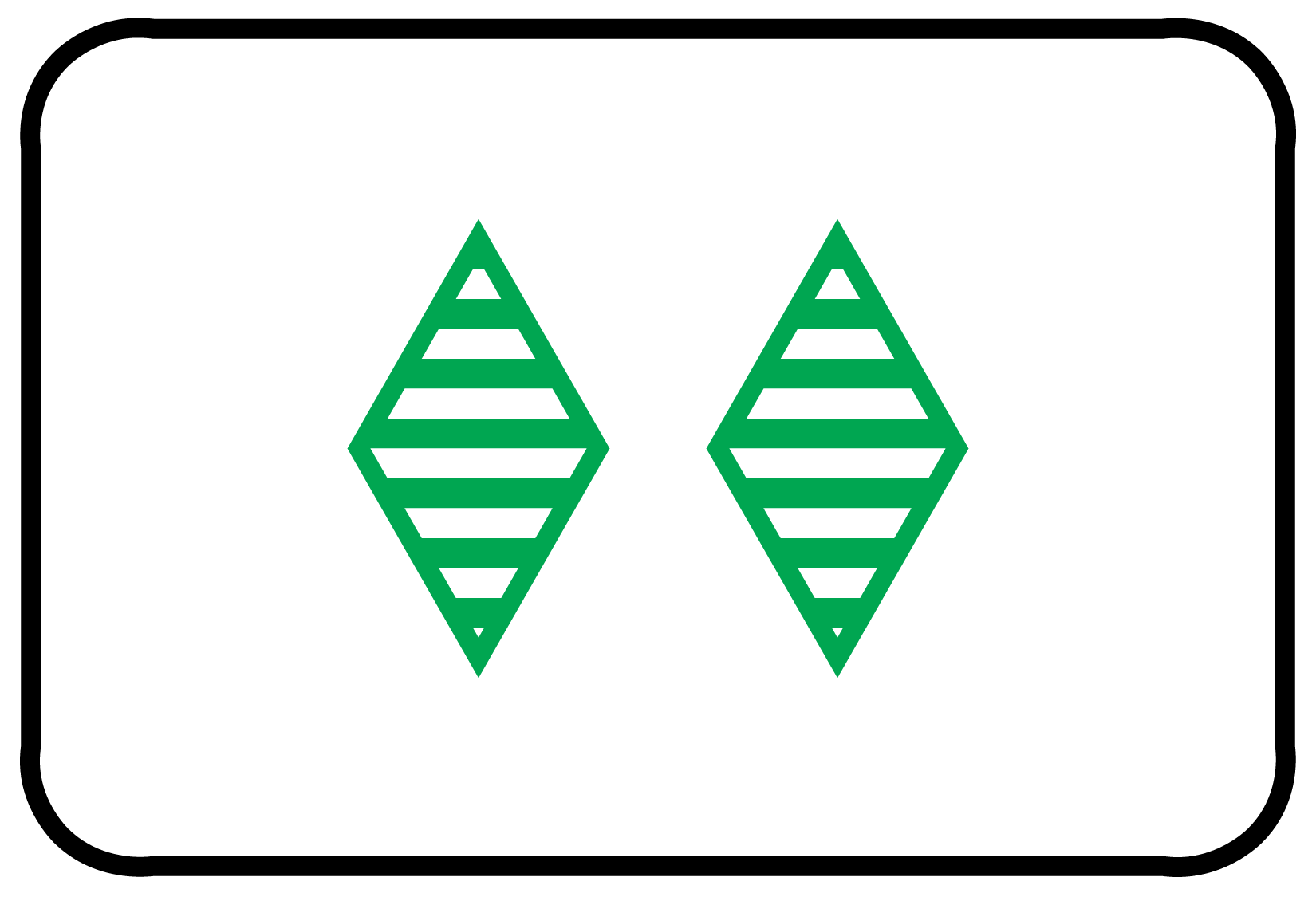} & \includegraphics[width = 1.25cm]{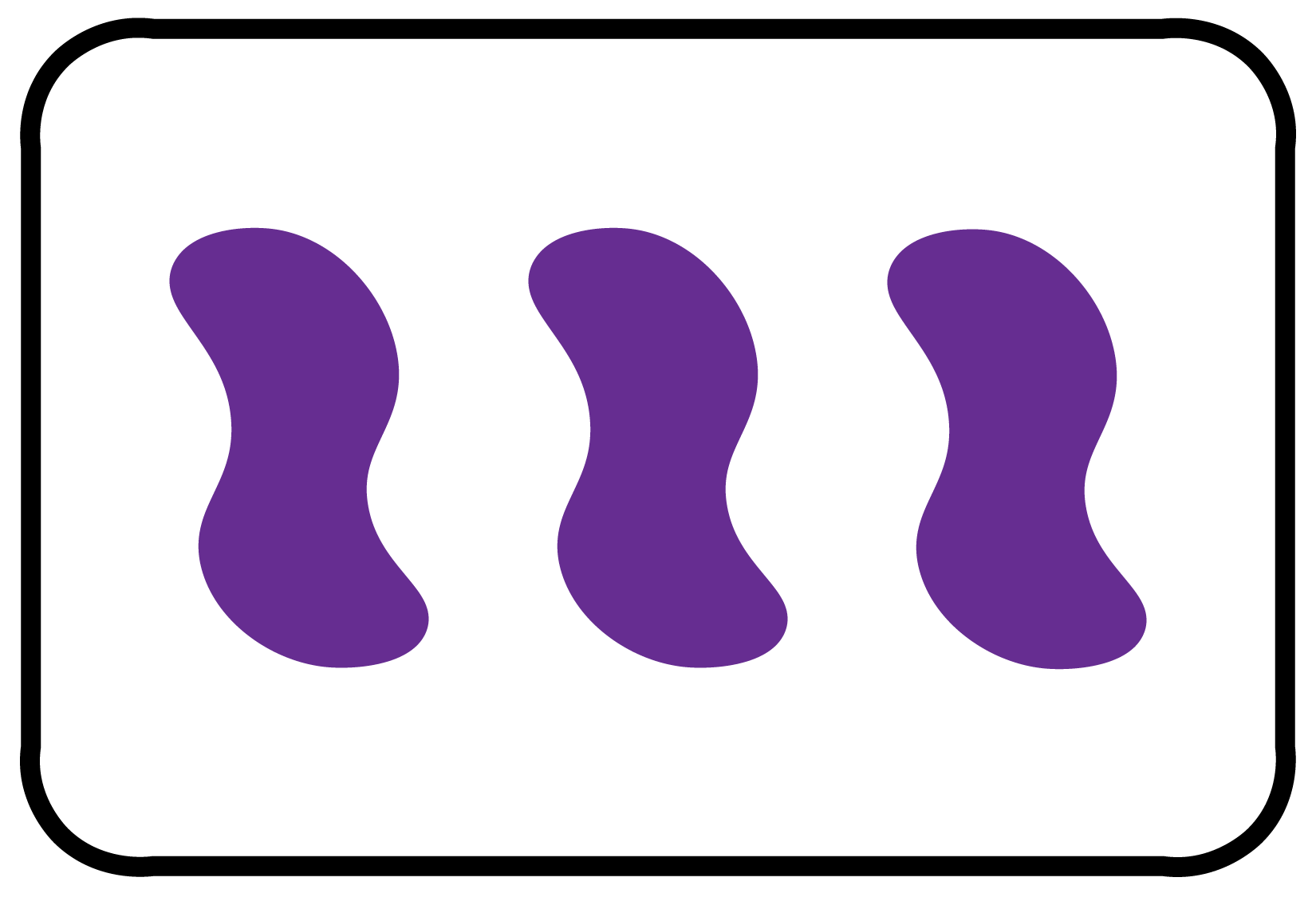}\\
        \end{tabular}
\caption{A set of order 4.}
\label{fig:set4}
   \end{minipage}
\end{figure}

Figure~\ref{fig:set1} shows a set of order 1, while Figure~\ref{fig:set4} shows a set of order 4.

In this paper, we are not interested in the game's rules, which can be easily found on the Internet. Suffice it to say, that the player who is the fastest to find sets among given cards wins. 

This game is very mathematical and provides a fun introduction to combinatorics, number theory, and linear algebra. One can read about this game in a wonderful book ``The Joy of SET'' published in 2016 \cite{MGGG}.

We can assign a value 0, 1, and 2 to different shapes, numbers, colors, and shadings. Then each card is represented as a set of four digits. In our internal numbering, we used the digits to represent features in the order we mentioned: shapes, numbers, colors, and shadings. Also, we assigned 0/1/2 for oval/diamond/squiggle shape; 0/1/2 for number 1/2/3; 0/1/2 for red/green/purple for color; 0/1/2 for empty/striped/filled shading. For example, the string $0000$ represents a card with one empty red oval.

We can view each digit as a coordinate. Our set of cards becomes a four-dimensional vector space. As we have only values 0, 1, 2, this four-dimensional vector space has values in a field of three elements $\mathbb{F}_3$ and can be denoted as $\mathbb{F}_3^4$.

What is the value of adding linear algebra to some cards? The beauty of the vector space is that we can add vectors. In terms of cards, adding two cards means adding each coordinate modulo 3. Now we have a way fancier definition of a set. Three cards form a set if and only if, their sum is 0. 

We can check that three numbers from the set 0, 1, and 2 sum to zero modulo three if and only if they are all the same or all different. That means three cards sum up to zero, if and only if, their value in each coordinate is either the same or different. This was the definition of a set from the start.

In this paper, we are not interested in the game itself, but rather in the magic SET squares.

\section{Magic squares and magic SET squares}

\begin{figure}[htp!]
        \begin{center}
            \includegraphics[scale=0.4]{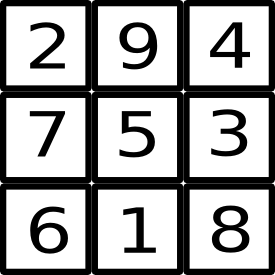}
        \end{center}
\caption{3 by 3 magic square.}
\label{fig:ms}
\end{figure}

Magic squares are squares made out of numbers so that each row, column, and diagonal sum up to the same number. Figure~\ref{fig:ms} shows an example of a 3 by 3 magic square using digits 1 through 9. The common sum in this square is 15.

\begin{figure}[ht!]
\centering
    \setsquare{0000}{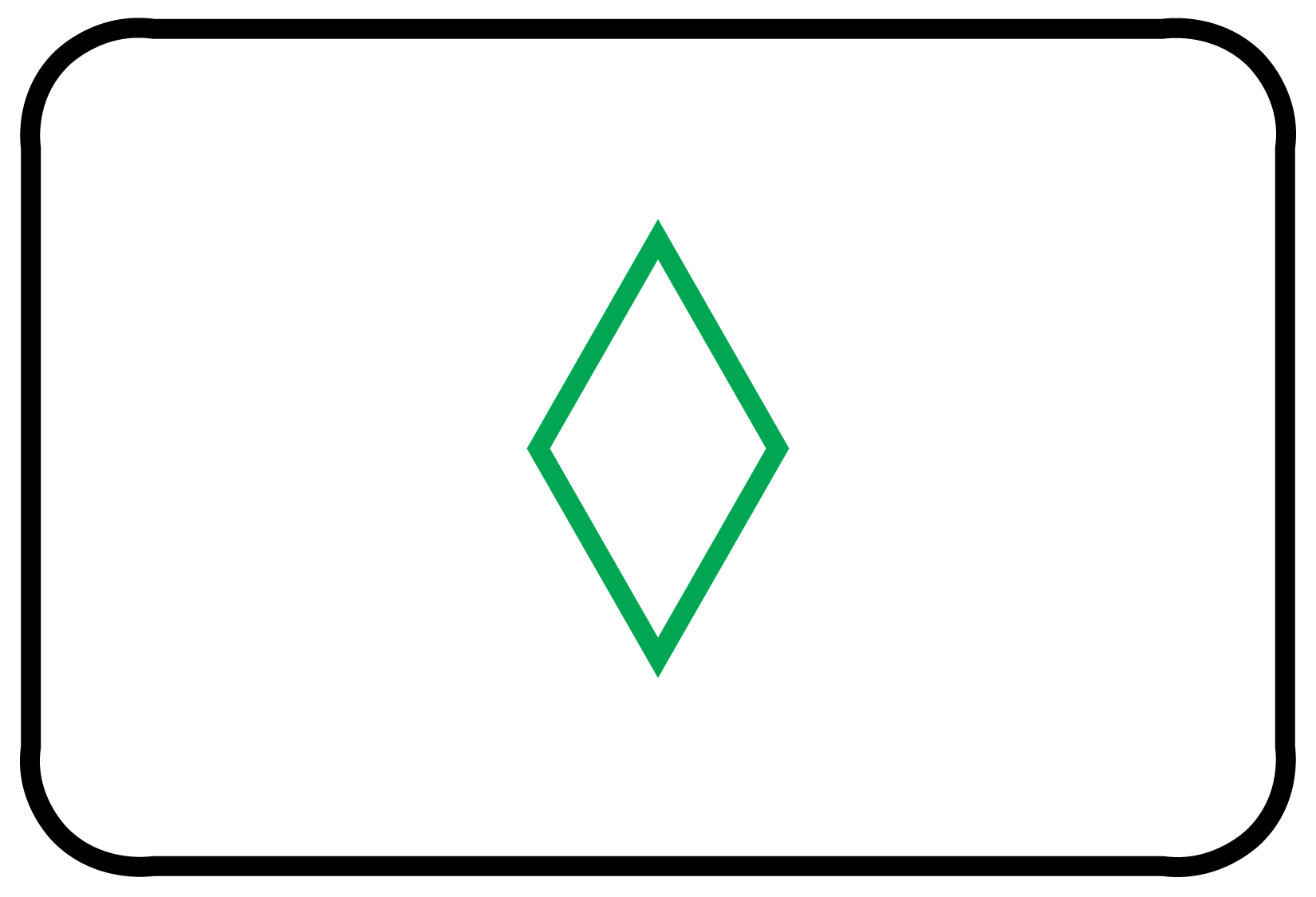}{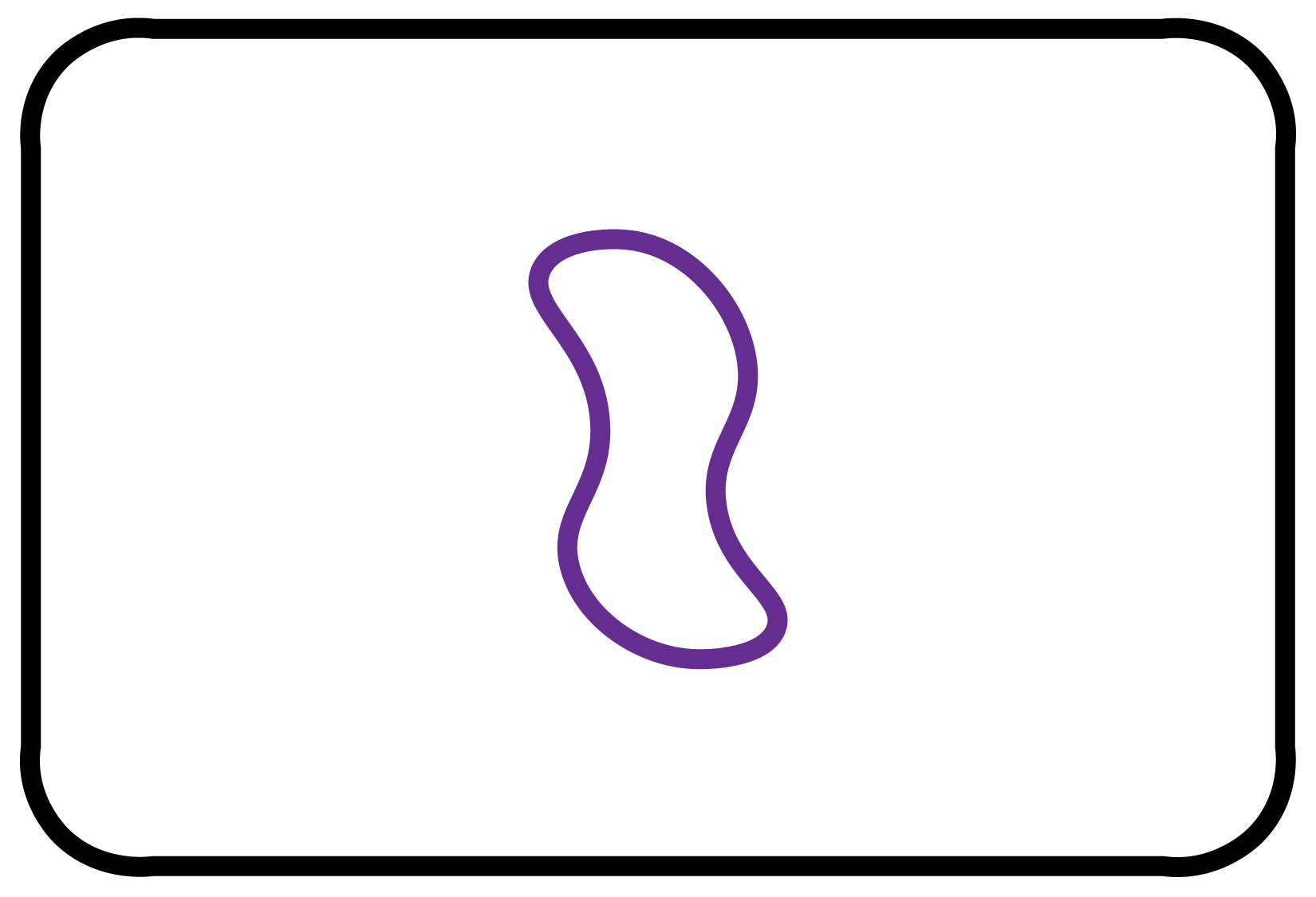}{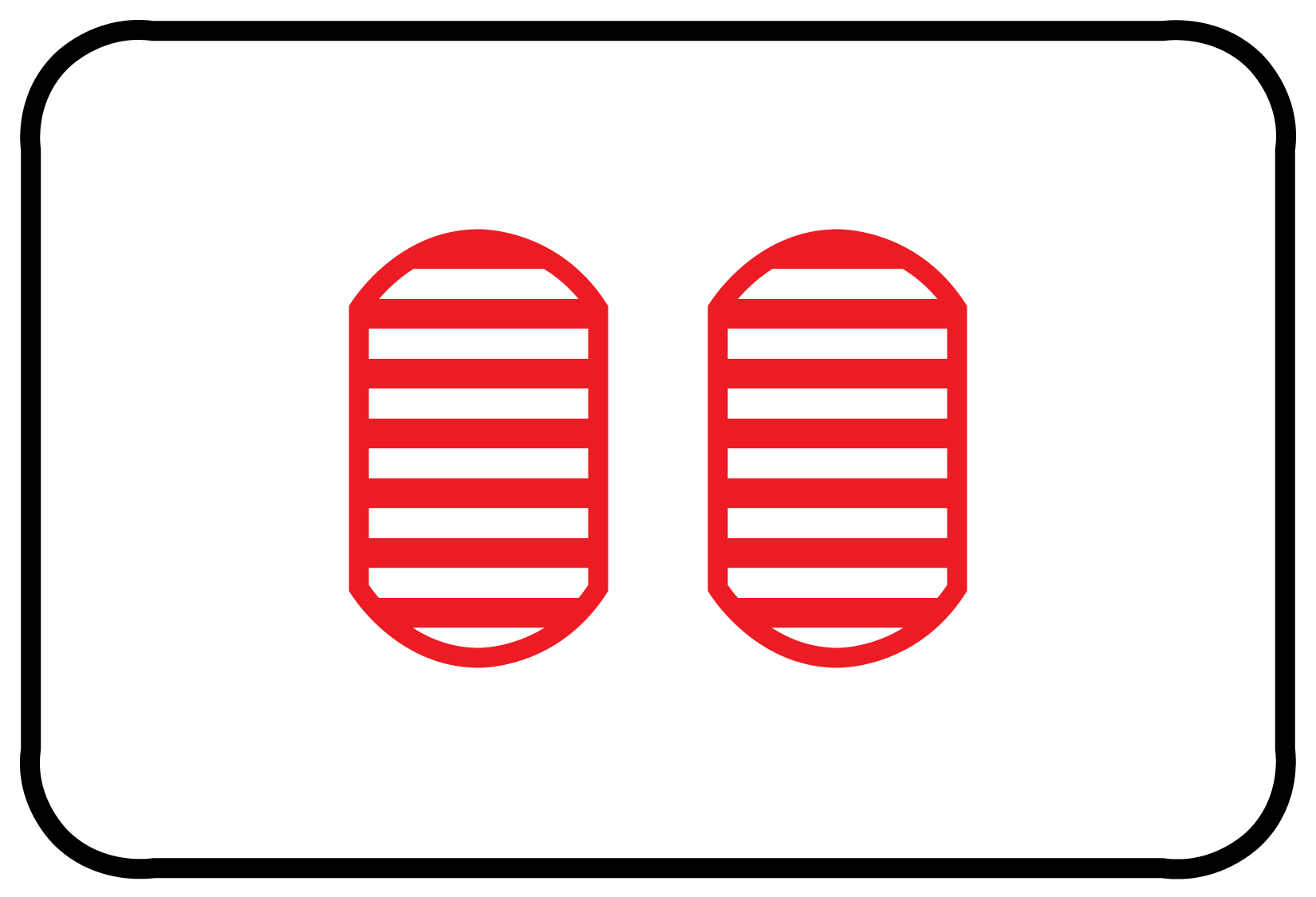}{1111}{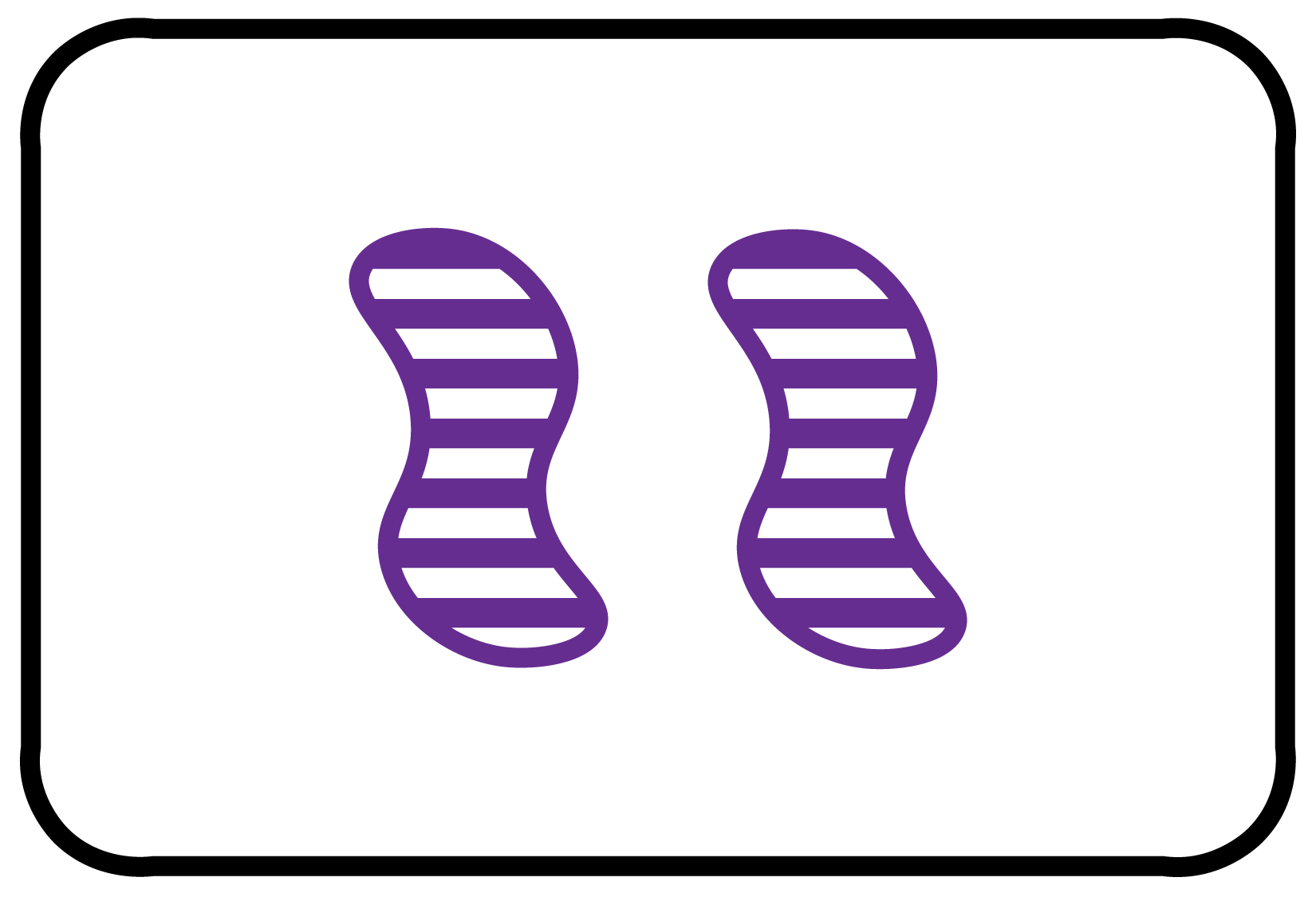}{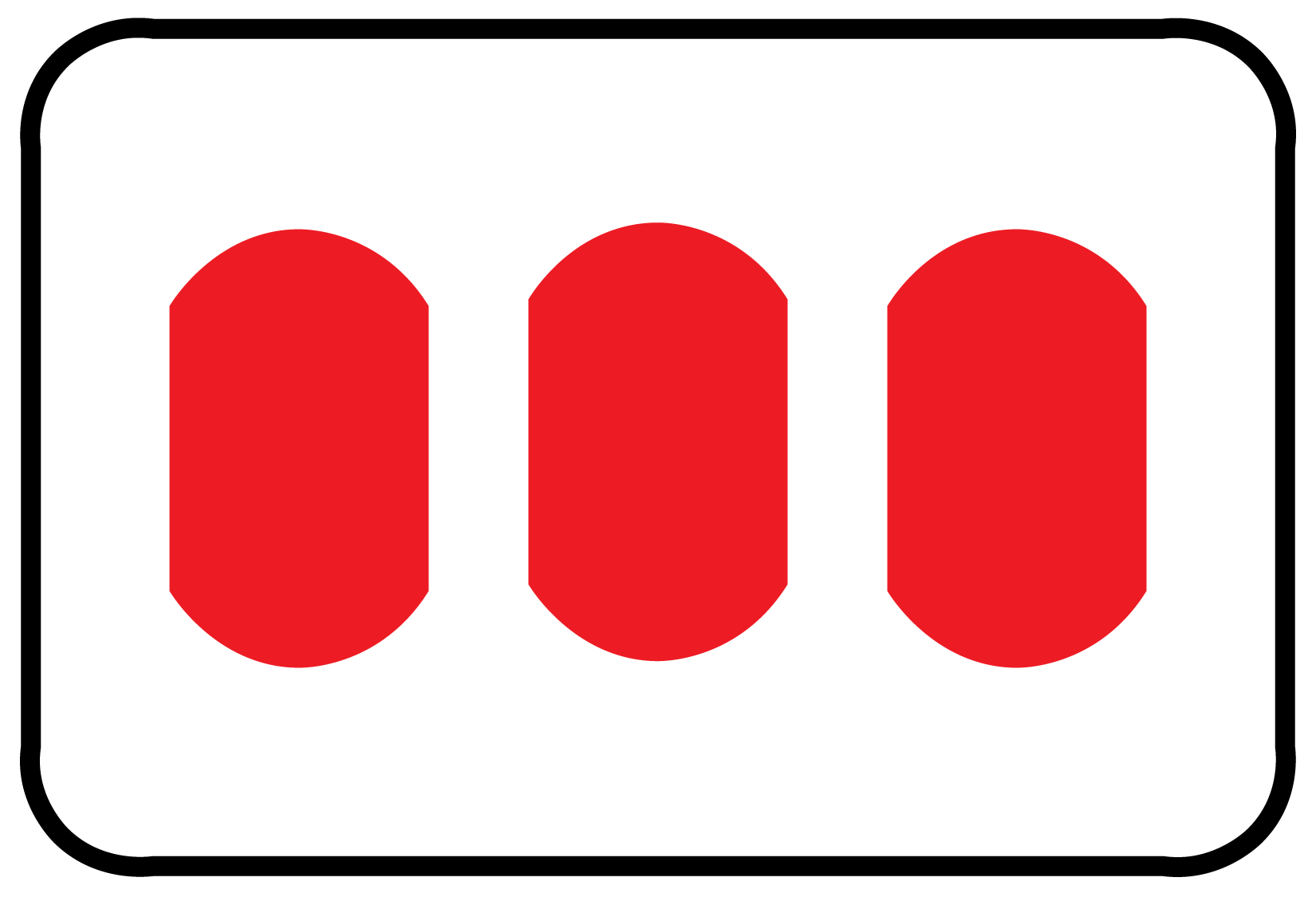}{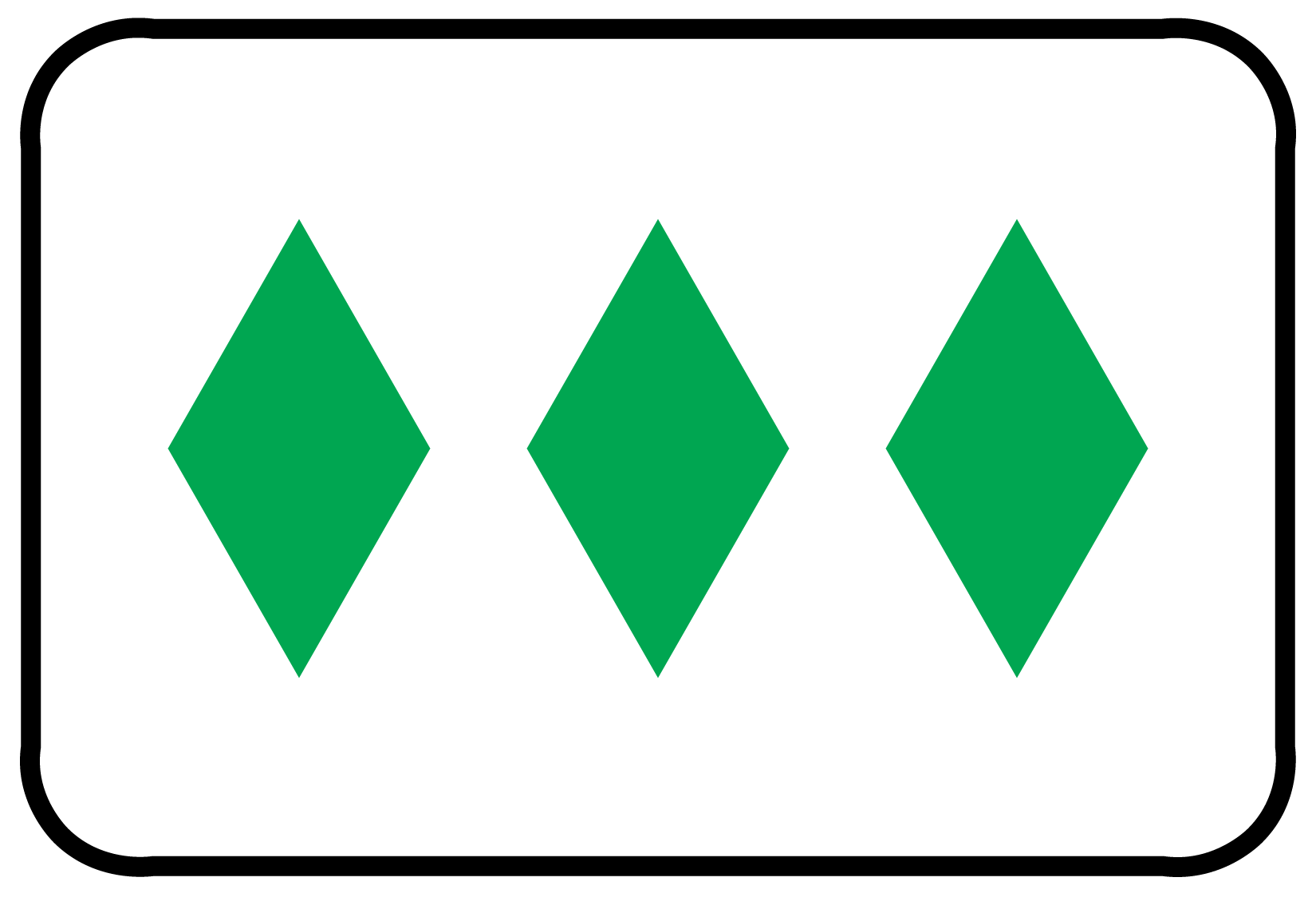}{2222}
\caption{A magic SET square.}
\label{fig:mssExample}
\end{figure}

We define a \textit{magic SET square} as a 3 by 3 table of SET cards, so that each row, column, and diagonal forms a set. Figure~\ref{fig:mssExample} shows an example of a magic SET square.

A magic SET square can be built starting with three cards, $a$, $b$, and $c$, that do not form a set. We place $a$ in the bottom left corner, $b$ in the bottom right corner and $c$ in the center of the square. We can calculate all the other cards in this square uniquely. The result is in Table~\ref{table:mss}.

\begin{table}[htp!]
\begin{center}
  \begin{tabular}{| c | c | c |}
    \hline
    $-b-c$ & $a+b-c$ & $-a-c$ \\ \hline
    $-a +b+c$ & $c$ & $a-b+c$ \\ \hline
    $a$ & $-a-b$ & $b$ \\
    \hline
  \end{tabular}
\end{center}
\caption{A magic SET square calculated from three cards $a$, $b$, and $c$}
\label{table:mss}
\end{table}

The reader can check that indeed every row, column, and diagonal forms a set. For example, let us see that the cards in the first row form a set. If we sum them up, we get $-3c$ which is zero in our addition.

\begin{figure}[ht!]
\centering
               \subfigure[Columns]{\includegraphics[scale=0.31]{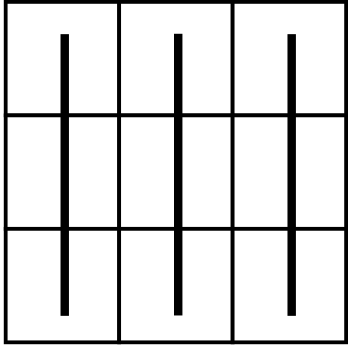}}\hfill
               \subfigure[Rows]{\includegraphics[scale=0.31]{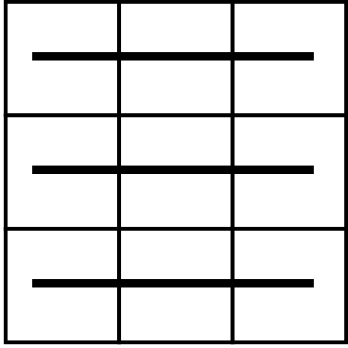}}\hfill
               \subfigure[Diagonals]{\includegraphics[scale=0.31]{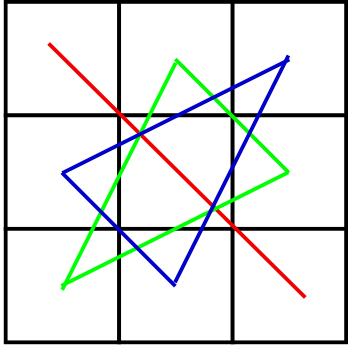}}\hfill
               \subfigure[Anti-diagonals]{\includegraphics[scale=0.31]{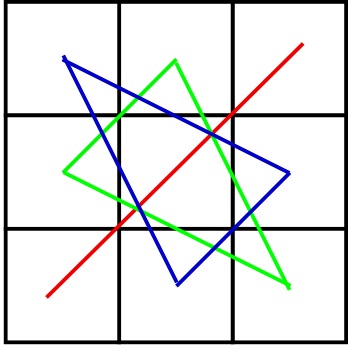}}
\caption{Triplets of sets in a magic SET square.}
\label{fig:triplets}
\end{figure}

If you play with magic SET squares you can notice that you can always find some more sets. The total number of sets in a magic SET square is always 12: three sets that form rows, three sets that form columns, three sets that form a diagonal and two broken diagonals, and three sets that form an anti-diagonal and two broken anti-diagonals. We call each group a \textit{triplet}. We have row, column, diagonal and anti-diagonal triplets. Figure~\ref{fig:triplets} shows the positions of triplets.

Similar to a set, we say that the \textit{order} of a magic SET square is the number of features that are different in the square. For example, the square in Figure~\ref{fig:mssExample} is of order 4 as we have cards of different colors, shapes, shadings, and numbers. 

The reader might notice that in our example of a magic SET square all the sets in one triplet have the same order. Is this always the case? Yes, it is always the case. 

In each triplet the sets are parallel. Two sets $A$ and $A'$ are \textit{parallel} if they differ by a vector. That is $A' = A+d$. For example the bottom row in Table~\ref{table:mss} is vector $\{a,-a-b,b\}$, while the middle row is $\{-a+b+c,c,a-b+c\}$. Subtracting the middle row from the bottom row we get $\{2a-b-c,-a-b-c,-a +2b -c\} = \{-a-b-c,-a-b-c,-a -b -c\}$.

Parallel sets have the same type for any feature. That is, all of them have either the same value or a different value for the same feature. The important conclusion is that parallel sets are of the same order. It follows that each triplet consists of three sets of the same order. You can read more about this in \cite{MGGG}.

We would like to differentiate squares of the same order. For this purpose, we introduce the notion of diversity. The \textit{diversity} of a magic SET square is the average order of all its sets. In the example in Figure~\ref{fig:mssExample} the rows are order 2: they differ in color and shape. The columns are order 2: they differ in number and shading. The diagonals and anti-diagonals are order 4. Thus the diversity of this square is 3.

Rows and columns are perceived differently. Thus, we introduce the row-column diversity, which we call \textit{rc-diversity} for short. The \textit{rc-diversity} of a magic SET square is the average order of the sets that form rows and columns of the square. In Figure~\ref{fig:mssExample} the rc-diversity is 2.

\section{The number of sets}

The number of sets in the deck is 1080. Now, we count sets of order $k$. First, we count ordered sets. There are $\binom{4}{k}$ ways to choose features for which the cards are different in this set. There are 3 ways to choose a value for each feature that is the same, resulting in $3^{4-k}$ total ways. For each feature that is different, there are 3 ways to choose the value for the first card and 2 ways for the second card. Thus, we need to multiply by $6^k$. The number of ordered sets is $\binom{4}{k}3^{4-k}6^k$. Then, for unordered sets, we can shuffle the cards in the set in 6 ways. The final answer is 
\[\binom{4}{k}3^{4-k}6^{k-1}.\]

Thus we have $4\cdot 3^3 \cdot 6^0 = 108$ sets of order 1. Also, $6\cdot 3^2 \cdot 6^1 = 324$ sets of order 2. Continuing, we have $4\cdot 3^1 \cdot 6^2 = 432$ sets of order 3. And, finally, we have $1\cdot 3^0 \cdot 6^3 = 216$ sets of order 4.

We put these numbers into Table~\ref{table:nsets}.

\begin{table}[htp!]
\begin{center}
\begin{tabular}{ |c|c|} 
 \hline
 Order & Number of sets \\ 
	\hline
 1 & 108 \\ 
 2 & 324 \\ 
 3 & 432  \\ 
 4 & 216  \\ 
\hline
 total & 1080 \\
 \hline
\end{tabular}
\end{center}
\caption{Number of sets of different orders.}
\label{table:nsets}
\end{table}

\section{The number of magic SET squares by order}

We can calculate the total number of magic SET squares as follows. We can pick any two cards for the bottom left and right corners in $81 \cdot 80$ ways. Then there are 78 ways for the center card as it shouldn't form a set with the bottom two cards. The total is: 505440.

The magic SET squares of order 1 do not exist. Indeed, if we fix the values for three features, then there are only three cards in a deck with these values. You can't form a 3 by 3 square out of three cards! Thus the magic SET square can only be of orders 2, 3, and 4.

Now we calculate the total number of squares of order 2. There are 2 characteristics that are the same. That means there are $\binom{4}{2}=6$ ways to pick these characteristics and nine ways to pick specific values for them. That means there are $6 \cdot 9 = 54$ different ways to pick the cards for the set. After we picked the cards, we can choose the bottom left corner in 9 ways, the bottom right corner in 8 ways, and the center in 6 ways. The total is $54 \cdot 9 \cdot 8 \cdot 6 = 23328$.

Now we calculate the number of squares of order at most 3. First, we choose the feature that is the same, then the value for this feature. There are 12 ways to do that. Now we have 27 cards to choose from for our square. We can choose the bottom left in 27 ways, the bottom right in 26 ways, and the center 24 ways. Recall that in our 27 cards any 2 cards can be completed into a set, and we need the card for the center that doesn't form the set with the bottom row. The result is 202176.

To get to the number of squares of order 3, we need to subtract the number of squares of order 2 from the number above. We counted each square of order 2 twice. The square of order 2 has two fixed values for 2 different features. We counted the square each time we fixed one feature. The final result for the number of squares of order 3 is 155520.

To find the number of squares of order 4 we take the total number of squares and subtract the number of squares of order 2 and 3. This is $505440 - 23328 -155520 = 326592$. We summarize the results in Table~\ref{table:nsquares}.

\begin{table}[htp!]
\begin{center}
\begin{tabular}{ |c|c|} 
 \hline
 Order & Number of squares \\ 
	\hline
 2 & 23328 \\ 
 3 & 155520  \\ 
 4 & 326592  \\ 
\hline
 total & 505440 \\
 \hline
\end{tabular}
\end{center}
\caption{Number of magic SET squares of different orders.}
\label{table:nsquares}
\end{table}

\section{Magic SET squares with the same cards}

Given a magic SET square, we can rearrange the cards to make another magic square. How many different squares can we make? First, we notice that for any two cards in a square, the card that completes them to a set is also in the square. That means we can pick any two cards for the bottom corners. There are $9 \times 8=72$ ways to do that. Then we can choose another card that doesn't form a set with bottom corners for the center. There are 6 ways to do that. So the number of different squares using the same cards is $72 \times 6=432$. Another way to calculate this is to choose a set for the bottom row in 12 ways, then shuffle the cards in the set in 6 ways, and then choose the center card out of the leftover 6 cards. We get the same total, as expected.

We want to discuss possible rearrangements of the cards in a square in more geometric terms.

As we mentioned before, sets in the same triplet are parallel to each other. In addition to that, one might notice that two sets from different triplets always share a card. That means they can't be parallel. It follows that when we rearrange the cards the triplets of sets stay together in a triplet. Now we want to show how we can shuffle triplets.

We can reflect a square along the main diagonal or the main anti-diagonal. This way rows and columns swap, but diagonal and anti-diagonal triplets stay in place. 

We can reflect a square along the middle column or the middle row. This way rows and columns stay in place and diagonals and anti-diagonals swap. 

We want to show that we can shuffle triplets in any possible way. 
We know that a permutation group can be generated by adjacent transpositions. We already showed how to swap rows and columns. Separately, we can swap diagonals and anti-diagonals. What is left to see is a transformation that swaps columns and diagonals.

We start with our square represented in Table~\ref{fig:trans1}.

\begin{table}[ht!]
\begin{center}
 \begin{tabular}{||c c c||} 
 \hline
 A & B & C  \\ 
 \hline
 A+X & B+X & C+X \\
 \hline
 A+2X & B+2X & C+2X \\
 \hline
\end{tabular}
\end{center}
\caption{Starting square.}
\label{fig:trans1}
\end{table}

Then we reflect it with respect to the middle column and get a square in Table~\ref{fig:trans2}.
\begin{table}[ht!]
\begin{center}
 \begin{tabular}{||c c c||} 
 \hline
 C & B & A  \\ 
 \hline
 C+X & B+X & A+X \\
 \hline
 C+2X & B+2X & A+2X \\
 \hline
\end{tabular}
\end{center}
\caption{After the reflection with respect to the middle column.}
\label{fig:trans2}
\end{table}

Then we shift the second and the third row to the right to make a parallelogram as seen in Table~\ref{fig:trans3}.
\begin{table}[ht!]
\begin{center}
 \begin{tabular}{||c c c c c||} 
 \hline
 C & B & A & & \\ 
 \hline
  & C+X & B+X & A+X  &\\
 \hline
  & & C+2X & B+2X & A+2X \\
 \hline
\end{tabular}
\end{center}
\caption{A shift.}
\label{fig:trans3}
\end{table}
and wrap it around as shown in Table~\ref{fig:trans4}.
\begin{table}[ht!]
\begin{center}
 \begin{tabular}{||c c c||} 
 \hline
 C & B & A\\ 
 \hline
  A+X & C+X & B+X \\
 \hline
  B+2X & A+2X & C+2X \\
 \hline
\end{tabular}
\end{center}
\caption{A wrap around.}
\label{fig:trans4}
\end{table}

In the resulting square row and anti-diagonal triplets stay in place, and column and diagonal triplets are swapped.
Therefore, since all adjacent transformations are possible, all permutations of the triplets are possible.

If we divide the total number of squares that we can make with the same cards by 24, the number of permutations of triplets, we get 18. This is the number of squares that we can make with the same cards keeping triplets in place.

One way to describe these 18 squares is the following. We can take a square and cycle rows or columns. This creates 9 squares out of one. We can also rotate a square 180 degrees. We need to mention that the rotation can not be achieved by cycling. Indeed, when we cycle rows or columns, the cards in rows or columns undergo an even permutation. During the rotation, the bottom row cards become the top row cards in the reverse order, which is an odd permutation of cards in each row or column. Thus, all possible transformations of a square that keep triplets in place can be described as a combination of row/column cycling and a 180-degree rotation of the square.

Another way to see this is the following. We can shuffle rows and columns in 6 ways each. Shuffling the bottom row and the left column defines a unique magic square. Thus, we create 36 squares where row and column triplets stay in place. Meanwhile, diagonal and anti-diagonals either stay in place or swap. That means in half of these squares all triplets stay in place.

\section{Our classification}

Consider the orders of triplets of a magic SET square $(r,c,d,a)$, where $r/c/d/a$ represent the orders of rows/columns/diagonals/anti-diagonals correspondingly. We denote the number of such squares as $N(r,c,d,a)$. We showed that there are transformations that can shuffle features. Suppose $s$ is such a transformation. As such shuffles are one-to-one, we have $N(s(r),s(c),s(d),s(a)) = N(r,c,d,a)$.

We consider magic SET squares with the same set of orders to be of the same type. We also consider two magic SET squares of the same type if one can be moved to another by a geometric transformation of the square, that is, by a rotation or reflection.

We consider two magic SET squares with triplet orders $(r_1,c_1,d_1,a_1)$ and $(r_2,c_2,d_2,a_2)$ to be of the same type if the sets $\{r_1,c_1\}$ and $\{r_2,c_2\}$ are the same as well as the sets $\{d_1,a_1\}$ and $\{d_2,a_2\}$ are the same. The reason is that we can swap row and column triplets by reflecting the square with respect to a diagonal or an anti-diagonal. This transformation keeps the diagonal and anti-diagonal triplets in place. We can swap the diagonal anti-diagonal triples with the reflection along the middle row or the middle column. This reflection keeps row and column triplets in place.

Thus we can represent the type of a square by four numbers $(x$-$y;z$-$w)$, where $x \leq y$ and $z \leq w$. The numbers $x$ and $y$ represent orders of rows and columns and it doesn't matter in which order we use them. The numbers $z$ and $w$ represent orders of diagonals and anti-diagonals and it doesn't matter in which order we use them. We denote the number of squares of type $(x$-$y;z$-$w)$ as $N(x$-$y;z$-$w)$.

What is the connection between the numbers $N(x,y,z,w)$ and $N(x$-$y;z$-$w)$? The answer: it depends. If $x=y$ and $z=w$, then both numbers count the same squares and we have $N(x,y,z,w) = N(x$-$y;z$-$w)$. If only one equality is true, that is, either $x=y$ or $z=w$, then $2N(x,y,z,w) = N(x$-$y;z$-$w)$. Indeed, suppose $x=y$ and $z \neq w$, then $N(x$-$y;z$-$w) = N(x,y,z,w) + N(x,y,w,z) = 2N(x,y,z,w)$.  If both equalities are false, that is, $x \neq y$ and $z \neq w$, then $4N(x,y,z,w) = N(x$-$y;z$-$w)$. Indeed, $N(x$-$y;z$-$w) = N(x,y,z,w) + N(y,x,z,w) + N(x,y,w,z) +N(y,x,w,z)$.

There is also symmetry between $(x$-$y;z$-$w)$ and $(z$-$w;x$-$y)$. We already showed how to swap columns and diagonals. In the same way, we can also swap rows and anti-diagonals. If we do both swaps, we have swapped the columns and rows with the diagonals and anti-diagonals correspondingly. Thus the number of magic SET squares of type $(x$-$y;z$-$w)$ is the same as the  number of squares of type $(z$-$w;x$-$y)$.

\section{Classification of magic SET squares of order 2}

We start by studying magic SET squares that have 9 cards with two features the same in all of them. We call such squares order-2 squares. This set of cards is uniquely defined by the values of these two features  that are the same. Looking at such nine cards one might see that the orders or triplets have to be $(1,1,2,2)$ up to a permutation.

All the squares we can make out of these cards have diversity 1.5, and there are three types of such squares: $(1-1;2-2)$, $(1-2;1-2)$, and $(2-2;1-1)$. A magic SET square of type (1-1;2-2) is in Figure~\ref{fig:(1-1;2-2)}, a magic SET square of type (2-2;1-1) is in Figure~\ref{fig:(2-2;1-1)}, and a magic SET square of type (1-2;1-2) is in Figure~\ref{fig:(1-2;1-2)}. We calculate the number of squares of each type in this section.

\begin{figure}[!ht]
   \begin{minipage}{0.3\textwidth}
\setsquarewithlabels{(1-1;2-2)}{0000}{1000}{2000}{0100}{1100}{2100}{0200}{1200}{2200}
   \end{minipage}\hfill
   \begin{minipage}{0.3\textwidth}
\setsquarewithlabels{(2-2;1-1)}{0000}{0011}{0022}{0021}{0002}{0010}{0012}{0020}{0001}
   \end{minipage}\hfill
   \begin{minipage}{0.3\textwidth}
\setsquarewithlabels{(1-2;1-2)}{2012}{2210}{2111}{2112}{2010}{2211}{2212}{2110}{2011}
   \end{minipage}
\end{figure}

We start by calculating the number of squares of type (1-1;2-2). Consider a square where each row and column are first order, that is rc-diversity is 1. We want to count the number of such squares. We choose the card in the left bottom corner in 81 ways. The card in the right bottom corner has to differ only in one feature, so there are 8 possibilities: 4 ways to choose the feature to be different and 2 ways to choose the value for this feature. Now we want to choose the rest of the square. For this, we need to choose the feature that is different in the column in 3 ways. This will uniquely define the cards in the columns. The only transformation that we can still do is to swap the top and middle rows. So there are 
\[81 \cdot 8 \cdot 6 = 3888\]
different magic SET squares of this type. 

Here is another calculation. There are 108 different sets of order 1. We can choose one of such sets. Then we can shuffle the cards in this set in 6 ways. So there are $108 \cdot 6 =648$ different ways to arrange the bottom row. After that, we can pick the center card so that it is the first-order set together with the middle card in the bottom row. There are 3 ways to choose a new feature in which this card is different, and 2 ways to choose the value of this card. Thus, the number of such sets is $648 \cdot 6 = 3888$.

Here is yet another calculation. We know that there are 23328 magic SET squares of order 2. We also know that there are 432 ways to arrange a square with the same cards. That means there are  $23328/432 = 54$ different ways to pick a set of cards for the magic SET square of order 2. Suppose we have the cards for the second-order square. Among these cards, we have 6 first-order sets. Suppose we pick one of the sets for the bottom row and then shuffle the cards for a particular order in the bottom row. Thus, there are 36 ways to pick the bottom row. Now we pick the card in the top left corner. The only condition is that it has to form a set of the first order with the bottom left card. There are 2 cards like this.  Thus there are $54 \cdot 36 \cdot 2 = 3888$ squares. 

The number of squares of type (2-2;1-1) is the same as the number of squares of (1-1;2-2). It is also the same as the number of squares of type $(1,1,2,2)$. The number of squares of type (1-2;1-2) is $N(1-2;1-2)$, which we calculate as:
\[N(1-2;1-2) = 4N(1,2,1,2) =4N(1,1,2,2)=4N(1-1;2-2) = 4 \cdot 3888 = 15552.\]

Magic squares of order 2 all have diversity 1.5. But there are three types of such squares and they all have different rc-diversity. The number of the squares of each type are summarized in Table~\ref{table:2bytype}.

\begin{table}[ht!]
\begin{center}
\begin{tabular}{ |c|c|c|} 
 \hline
Type & rc-diversity & \#squares\\
\hline
(1-1;2-2) & 1 & 3888  \\ 
(1-2;1-2) & 1.5 & 15552  \\ 
(2-2;1-1) & 2 & 3888  \\ 
\hline
\end{tabular}
\caption{Number of squares of order 2 by type.}
\label{table:2bytype}
\end{center}
\end{table}

The total number of magic squares of order 2 matches the number we calculated before: 23328.

\section{Order of a magic SET square and its diversity}

We noticed that all magic SET squares of order 2 have the same diversity. Is there a reason for that? It is often a great idea to look at each feature separately. Consider a magic SET square. What can we say about how one feature can be distributed? Obviously, all the cards can be the same in this feature. If not, what can we say about the triplets? Suppose nine cards in the square have different shapes. It means there are three diamonds, three squiggles, and three ovals in the square. It means in terms of this feature only one triplet is the same and three triplets are different. This feature contributes $\frac{3}{4}$ to the diversity of the square. 

Thus, each feature in which not all cards in the square are the same contributes $\frac{3}{4}$ to the diversity of the square. Remember that a square of order $k$ has $k$ features that are different. That means such a square has diversity $\frac{3k}{4}$. Table~\ref{table:orderdiversity} shows the explicit value for the diversity of a square depending on the order.

\begin{table}[ht!]
\begin{center}
\begin{tabular}{ |c|c|} 
 \hline
Order & diversity \\
\hline
2 & 1.5 \\ 
3 & 2.25 \\ 
4 & 3  \\ 
\hline
\end{tabular}
\caption{Order and diversity.}
\label{table:orderdiversity}
\end{center}
\end{table}

Now we can go back and see how the orders of triplets must be distributed in each square. The order of a triplet in a square of order $k$ must be between 1 and $k$ inclusive.

So, in a square of order 2, the orders of  triplets must be either 1 or 2 and have to sum up to 6. The only possibility is 1, 1, 2, and 2.

Consider a square of order 3. Each triplet can have order 1, 2, or 3, and the sum of orders must be 9. There are two possibilities: $\{1,2,3,3\}$ and $\{2,2,2,3\}$. For squares of order 4, the sum must be 12. Thus we have the following options: $\{1,3,4,4\}$, $\{2,2,4,4\}$, $\{2,3,3,4\}$, and $\{3,3,3,3\}$.

\section{Magic SET squares of order 3}

We showed that for magic SET squares of order 3, the orders of triplets are either $\{1,2,3,3\}$ or $\{2,2,2,3\}$. 

For the sake of practice and example, we show directly that if a square of order 3 has a triplet of order 1, then the orders of triplets must be $\{1,2,3,3\}$. Without loss of generality, we can assume that rows have order one.

Suppose the columns have order 2. The features that are different in row/column triplets have to be complementary, or the square will not be of order 3. Suppose we take two cards from different rows and columns. They have to differ in all three features. That means the diagonals and anti-diagonals are of order 3.

Suppose the columns have order 3. Denote by  $a$ the feature that is different in both columns and rows, and by $b$ and $c$ the features that are different in the columns but the same in the rows. If we pick any two cards that are neither in the same column nor in the same row, they have to differ in both $b$ and $c$. But they might be either the same of different in feature $a$. That means the other two triplets must be of orders 2 and 3. We proved that if there is a set of order 1 in a square of order 3, then the orders of triplets must be 1, 2, 3, and 3. Looking for each feature separately and using the diversity to find possible distributions was a simpler way to do this.

Thus, if there is a set of order one, the orders of triplets must be 1, 2, 3, 3, and the ones with orders 1 and 2 must have to be different on complementary features. The diversity of such a square is 2.25. Such squares are of four types: (1-2;3-3), (1-3;2-3), (2-3;1-2) and (3-3;1-2). Examples of such squares are shown in Figure~\ref{fig:(1-2;3-3)}, Figure~\ref{fig:(3-3;1-2)}, Figure~\ref{fig:(1-3;2-3)}, and Figure~\ref{fig:(2-3;1-3)} correspondingly.

\begin{figure}[!ht]
   \begin{minipage}{0.48\textwidth}
\setsquarewithlabels{(1-2;3-3)}{1002}{1021}{1010}{1102}{1121}{1110}{1202}{1221}{1210}
   \end{minipage}\hfill
   \begin{minipage}{0.48\textwidth}
\setsquarewithlabels{(3-3;1-2)}{0001}{2221}{1111}{1211}{0101}{2021}{2121}{1011}{0201}
   \end{minipage}
\end{figure}

\begin{figure}[!ht]
   \begin{minipage}{0.48\textwidth}
\setsquarewithlabels{(1-3;2-3)}{1012}{2111}{0210}{2012}{0111}{1210}{0012}{1111}{2210}
   \end{minipage}\hfill
   \begin{minipage}{0.48\textwidth}
\setsquarewithlabels{(2-3;1-3)}{0101}{0202}{0000}{0222}{0020}{0121}{0010}{0111}{0212}
   \end{minipage}
\end{figure}

First, we calculate the number of magic SET squares of order 3 so that rows are of order 1 and columns are of order 2. There are 81 ways to choose the bottom left card. After that, there are 8 ways to choose the bottom right corner: 4 ways to choose the feature that changes in the set and 2 more ways to pick the card from the set for the corner. After that, we need to pick 2 features that differ for the column. They can't include the feature that is different for the row, as otherwise the square will be of order 2. So there are 3 ways to pick these features. And there are four more ways to pick the values of the features for the top left card. Hence there are 12 ways to pick the top left card, so we have: $N(1,2,3,3) = 81 \cdot 8 \cdot 12 = 7776$. From here we get 
\[N(1-2;3-3) = N(3-3;1-2)= 2N(1,2,3,3) = 15552\]
and 
\[N(1-3;2-3) = N(2-3;1-3)= 4N(1,2,3,3) = 31104.\]

Now we consider squares of order 3 that do not contain a set of order 1. We found that there are squares where the orders in the groups are 2, 2, 2, and 3. Such squares are of two types: (2-2;2-3) and (2-3;2-2). Examples of such squares are shown in Figure~\ref{fig:(2-2;2-3)} and Figure~\ref{fig:(2-3;2-2)} correspondingly.

\begin{figure}[!ht]
   \begin{minipage}{0.48\textwidth}
\setsquarewithlabels{(2-2;2-3)}{0202}{0012}{0122}{0101}{0211}{0021}{0000}{0110}{0220}
   \end{minipage}\hfill
   \begin{minipage}{0.48\textwidth}
\setsquarewithlabels{(2-3;2-2)}{1001}{2010}{0022}{0000}{1012}{2021}{2002}{0011}{1020}
   \end{minipage}
\end{figure}

First, we calculate the number of squares where rows are order 3, and columns are order 2. We can pick a set of order 3 for the bottom row in 432 ways, then we can shuffle the cards in the set in 6 ways. Then, we choose the card in the center of the square. It differs from the bottom center card in two features. We can choose these two features out of all three features available in 3 ways. After that, we can choose the specific values for these features in 4 ways. Thus we have $ 432 \cdot 6 \cdot 3 \cdot 4 = 31104$ of such squares. These squares include the squares such that the diagonals and anti-diagonals are of order 1 and 3. Thus we need to subtract the number of such squares, which is 15552. We get the answer $N(3,2,2,2) = 15552$. From here we have 
\[N(2-3;2-2) = N(2-2;2-3) = 2N(3,2,2,2) = 31104.\]

We summarize the results in Table~\ref{table:3bytype}.

\begin{table}[ht!]
\begin{center}
\begin{tabular}{ |c|c|c|} 
 \hline
Type & rc-diversity & \#squares\\
\hline
(1-2;3-3) & 1.5 & 15552  \\ 
(1-3;2-3) & 2 & 31104  \\ 
(2-3;1-3) & 2.5 & 31104  \\ 
(3-3;1-2) & 3 & 15552 \\ 
(2-2;2-3) & 2 & 31104 \\ 
(2-3;2-2) & 2.5 & 31104 \\ 
\hline
\end{tabular}
\caption{Number of squares of order 3 by type.}
\label{table:3bytype}
\end{center}
\end{table}

The total number of squares of order 3 that contain a set of order 1 is the sum of the first three rows: $2\cdot 15552 +2 \cdot 31104 = 93312$. The number of squares of order 3 that do not contain a set of order 1 is 62208. The total number of magic SET squares of order 3 is 155520.

In our calculations, we are sometimes interested in features that are the same or different for a given pair of triplets. For example, in a square with triplets of order 2, 2, 2, and 3, consider two triplets of order 2. As the whole square is of order 3, the features that are different in these two triplets must share one feature in which they both are different. We have three triplets of order 2, any two triplets share a feature. This shared feature is different for every pair of triplets of order 2. We explain why this happens in the next section.

\section{Sharing features between triplets}

Now we want to discuss how different triplets share different features. 

Suppose one triplet has order $a$ and another triplet has order $b$, and the order of the square is $k$. If both triplets are the same in one feature, then the whole square will have the fixed value in this feature. That means the number of features that are different in at least one of the triplets is the order of the square. Thus, the number of features that are different in both triplets must be $a+b -k$. We call the corresponding features \textit{both-different features} for the pair of triplets.

For one feature that is different in a square, exactly three pairs of triplets are both-different. Moreover, if triplets $A$ and $B$ are both-different in this feature as well as triplets $A$ and $C$, then $B$ and $C$ are also both different and triplet $D$ has to be the same in this feature.

For example, consider a magic SET square of order 3 with orders of triplets $a$, $b$, $c$, and $d$ of orders 2, 2, 2, and 3 correspondingly. Then two triplets $a$ and $b$ of order 2 must be both-different in one feature, say $x$. But the set of order 3 has to be different in $x$ too. That means that the triplet $c$ is the same in feature $x$. We can have a similar argument for other pairs of triplets of order 2. It follows that three pairs of triplets of order 2 are both-different in three different features. Let us denote the three features that are different in this square of order 3 as $x$, $y$, and $z$. Then the following is the distribution of how the features are different in the triplets: $xy$, $xz$, $yz$, and $xyz$.

In Table~\ref{table:orderfeature} we show the distribution of features for all types of squares. This is done by using two principles: a) each feature is either the same in all triplets or is different in exactly three triplets out of four and b) for any two triplets the size of the union of all features that are different has to equal to the order of the square.

\begin{table}[ht!]
\begin{center}
\begin{tabular}{|c |c|c|} 
 \hline
Square order & Order dist. & Feature dist.\\
\hline
2 & 1, 1, 2, 2 & $x$, $y$, $xy$, $xy$ \\ 
3 & 1, 2, 3, 3 & $x$, $yz$, $xyz$, $xyz$ \\ 
3 & 2, 2, 2, 3 & $xy$, $xz$, $yz$, $xyz$ \\ 
4 & 1, 3, 4, 4 & $x$, $yzw$, $xyzw$, $xyzw$ \\ 
4 & 2, 2, 4, 4 & $xy$, $zw$, $xyzw$, $xyzw$ \\ 
4 & 2, 3, 3, 4 & $xy$, $xzw$, $yzw$, $xyzw$ \\ 
4 & 3, 3, 3, 3 & $xyz$, $xyw$, $xzw$, $yzw$ \\ 
\hline
\end{tabular}
\caption{Feature distribution.}
\label{table:orderfeature}
\end{center}
\end{table}

\section{Magic SET squares of order 4}

We have four cases for distribution of orders between all four triplets: $\{1,3,4,4\}$, $\{2,2,4,4\}$, $\{2,3,3,4\}$, and $\{3,3,3,3\}$. Now we give examples and calculate the number of sets of each type.

\subsection{Triplets have orders 1, 3, 4, 4}

Figure~\ref{fig:(1-3;4-4)}, Figure~\ref{fig:(4-4;1-3)}, Figure~\ref{fig:(1-4;3-4)}, and Figure~\ref{fig:(3-4;1-4)}   show examples of magic SET squares of type (1-3;4-4), (4-4;1-3), (1-4;3-4), and (3-4;1-4) correspondingly.

\begin{figure}[!ht]
   \begin{minipage}{0.48\textwidth}
\setsquarewithlabels{(1-3;4-4)}{2121}{0112}{1100}{2021}{0012}{1000}{2221}{0212}{1200}
   \end{minipage}\hfill
   \begin{minipage}{0.48\textwidth}
\setsquarewithlabels{(4-4;1-3)}{0011}{1120}{2202}{1220}{2002}{0111}{2102}{0211}{1020}
   \end{minipage}
\end{figure}

\begin{figure}[!ht]
   \begin{minipage}{0.48\textwidth}
\setsquarewithlabels{(1-4;3-4)}{0000}{1111}{2222}{0001}{1112}{2220}{0002}{1110}{2221}
   \end{minipage}\hfill
   \begin{minipage}{0.48\textwidth}
\setsquarewithlabels{(3-4;1-4)}{2022}{1111}{0200}{0000}{2122}{1211}{1011}{0100}{2222}
   \end{minipage}
\end{figure}

First, we calculate the number of squares where each row is order 1, and each column is order 3. There are 108 ways to choose a set of order 1 for the bottom row and 6 ways of ordering the cards in the set. To choose an order 3 column, the one common feature
in the column must be the feature that is different in the order 1 row. That means that there are $2 \cdot 2 \cdot 2 = 8$ ways
to choose a particular card in the column, as each of the 3 different features has 2
possibilities left for the value. Multiplying these all together, we get 5184.

This allows us to calculate the number of sets of each type. The results are in Table~\ref{table:1344}.

\begin{table}[ht!]
\begin{center}
\begin{tabular}{ |c|c|c|} 
 \hline
Type & rc-diversity & \#squares \\
\hline
(1-3;4-4) & 2 & 10368   \\ 
(1-4;3-4) & 2.5 & 20736  \\ 
(3-4;1-4) & 2.5 & 20736   \\ 
(4-4;1-3) & 3 & 10368   \\ 
\hline
total & --- & 62208  \\ 
\hline
\end{tabular}
\caption{Number of squares of different types for distribution 1, 3, 4, 4.}
\label{table:1344}
\end{center}
\end{table}

\subsection{Triplets have orders 2, 2, 4, 4}

Figure~\ref{fig:(2-2;4-4)}, Figure~\ref{fig:(4-4;2-2)}, and Figure~\ref{fig:(2-4;2-4)} show examples of magic SET squares of type (2-2;4-4), (4-4;2-2), and (2-4;2-4).

\begin{figure}[!ht]
   \begin{minipage}{0.3\textwidth}
\setsquarewithlabels{(2-2;4-4)}{1111}{0101}{2121}{1010}{0000}{2020}{1212}{0202}{2222}
   \end{minipage}\hfill
   \begin{minipage}{0.3\textwidth}
\setsquarewithlabels{(4-4;2-2)}{0202}{2020}{1111}{1010}{0101}{2222}{2121}{1212}{0000}
   \end{minipage}\hfill
   \begin{minipage}{0.3\textwidth}
\setsquarewithlabels{(2-4;2-4)}{1021}{0202}{2110}{1111}{0022}{2200}{1201}{0112}{2020}
   \end{minipage}
\end{figure}

To start off, we calculate the number of squares where the rows and columns are order 2. We can choose an order 2 set for the bottom row in 324 ways, and there are 6 ways to order the cards. The features in the order 2 column that are the same must be the features that are different in a row to make the square to be order 4. There are $2 \cdot 2 = 4$ cards that can go in a particular spot in the column, as there are 2 possibilities for values of each feature that are
different. Multiplying this all together, we get 7776.

This allows us to calculate the number of sets of each type. The results are in Table~\ref{table:2244}.

\begin{table}[ht!]
\begin{center}
\begin{tabular}{ |c|c|c|c|} 
 \hline
Type & rc-diversity & \#squares \\
\hline
(2-2;4-4) & 2 & 7776   \\ 
(2-4;2-4) & 2.5 & 31104   \\ 
(4-4;2-2) & 3 & 7776   \\ 
\hline
total & --- & 46656   \\ 
\hline
\end{tabular}
\caption{Number of squares of different types for distribution 2, 2, 4, 4.}
\label{table:2244}
\end{center}
\end{table}

\subsection{Triplets have orders 2, 3, 3, 4}

Figure~\ref{fig:(2-3;3-4)}, Figure~\ref{fig:(3-4;2-3)}, Figure~\ref{fig:(2-4;3-3)}, and Figure~\ref{fig:(3-3;2-4)} show examples of magic SET squares of type (2-3;3-4), (3-4;2-3), (2-4;3-3), and (3-3;2-4) correspondingly.

\begin{figure}[!ht]
   \begin{minipage}{0.48\textwidth}
\setsquarewithlabels{(2-3;3-4)}{0000}{0011}{0022}{1110}{1121}{1102}{2220}{2201}{2212}
   \end{minipage}\hfill
   \begin{minipage}{0.48\textwidth}
\setsquarewithlabels{(3-4;2-3)}{0022}{1110}{2201}{0111}{1202}{2020}{0200}{1021}{2112}
   \end{minipage}
\end{figure}

\begin{figure}[!ht]
   \begin{minipage}{0.48\textwidth}
\setsquarewithlabels{(2-4;3-3)}{0001}{0010}{0022}{1112}{1121}{1100}{2220}{2202}{2211}
   \end{minipage}\hfill
   \begin{minipage}{0.48\textwidth}
\setsquarewithlabels{(3-3;2-4)}{2210}{0012}{1111}{0000}{1102}{2201}{1120}{2222}{0021}
   \end{minipage}
\end{figure}

To start off, we calculate the number of such squares with rows of order 2 and columns of order 4. We can choose an order 2 set for the bottom row in 324 ways, and there are 6 ways to order the cards. The card in the center needs to have all features different from the center card in the bottom row. There are $2^4 = 16$ cards that can go in the center. Multiplying this all together, we get 31104. These squares include squares with orders (2,4,2,4) and (2,4,4,2). The number of such squares is a half of the number of squares of type (2-4;2-4). Thus, the total number of squares with orders (2,4,3,3) is $31104 - \frac{31104}{2} = 15552$. From here we get the summary in Table~\ref{table:2334},

\begin{table}[ht!]
\begin{center}
\begin{tabular}{ |c|c|c|} 
 \hline
Type & rc-diversity & \#squares \\
\hline
(2-3;3-4) & 2.5 & 62208   \\ 
(2-4;3-3) & 3 & 31104   \\ 
(3-3;2-4) & 3 & 31104  \\ 
(3-4;2-3) & 3.5 & 62208   \\ 
\hline
total & --- & 186624   \\ 
\hline
\end{tabular}
\end{center}
\caption{Number of squares of different types for distribution 2, 3, 3, 4.}
\label{table:2334}
\end{table}

\subsection{Triplets have orders 3, 3, 3, 3}

Figure~\ref{fig:(3-3;3-3)} shows an example of a magic SET square of type (3-3;3-3).

\begin{figure}[ht!]
\setsquarewithlabels{(3-3;3-3)}{1002}{2011}{0020}{2122}{0101}{1110}{0212}{1221}{2200}
\end{figure}

Suppose the bottom row has feature $x$ the same and other features different. Now we are choosing the center card of the square. It must be different in feature $y$ from the bottom row to make the square of order 4. Suppose it is the same with the bottom left card in feature $z$ that is different from $y$ and $x$. Then the center card is the same with the bottom right card in feature $w$ that is different from $x$, $y$, and $w$. There are 6 ways to choose the features for $y$, $z$, and $w$. After that, there are 2 ways to choose the value of the center card for feature $z$. The values for features $y$, $z$, and $w$ are uniquely defined. Overall, there are 12 ways to choose the center card, after the bottom row of order 3 is placed, for which there are $432 \cdot 6 = 2592$ ways. The total is $2592 \cdot 12 = 31104$, see Table~\ref{table:3333}.

\begin{table}[ht!]
\begin{center}
\begin{tabular}{ |c|c|c|} 
 \hline
Type & rc-diversity & \#squares \\
\hline
(3-3;3-3) & 3 & 31104   \\ 
\hline
total & --- & 31104   \\ 
\hline
\end{tabular}
\end{center}
\caption{Number of squares of different types for distribution 2, 3, 3, 4.}
\label{table:3333}
\end{table}

By summing the totals for each distribution, we get the expected total number of squares of order 4: 
\[62208 + 46656 + 186624 + 31104 = 326592.\]

\section{SET tic-tac-toe}

Everyone knows how to play tic-tac-toe. But what is its connection to magic squares? If you look at the magic square in Figure~\ref{fig:ms}, you can see that the straight lines in the square are in one-to-one correspondence with sets of three integers from 1 to 9 that sum to 15. That means that the game of tic-tac-toe is equivalent to the following game: Players start with a pool of nine cards that have the digits from 1 to 9 on them. They take turns taking cards from the pool. The first player who owns three cards that sum up to 15 wins. 

Now we can define SET tic-tac-toe. Clearly, it should be played using a magic SET square as the pool of cards. Similarly to the game of 15, the players take turns picking a card from the pool. The first player who has a SET wins. Equivalently, we can look at it as tic-tac-toe where broken diagonals and anti-diagonals are considered a line. 

Because broken diagonals and anti-diagonals are not valid in the game of 15, SET tic-tac-toe is not the same as the game of 15. Unlike regular tic-tac-toe, SET tic-tac-toe never ends in a draw. Indeed, if all the cards are drawn, the first player has five cards. But it is known that the largest number of cards in a magic SET square that do not contain a set is 4. Thus, the first player must have a set.

Now we show that the first player can always win after taking not more than 4 cards. Let us call the first player Alice, and the second player Bob. After each of them makes two turns, each player has two cards. If Alice can choose the card that completes her two cards to a set, then she wins on her third card. Suppose this is not the case. Then Alice, on her third turn, has to pick the card that completes the set for the second player, if this card is available. Otherwise, Alice can pick any card. Now Alice has three cards that do not form a set, and Bob can't complete a set on the next move, as Alice has that card.

Any two cards have a third card in the square that completes them to a set. As Alice has three cards not forming a set, there are three cards that can allow Alice to have a set, one per each pair. What is left to show is that Bob can't have all of these three cards.

Let Alice's cards be $a$, $b_1$, and $b_2$, where $a$ is the card that we know completes a set with two of Bob's cards. Let those be $c_1$ and $c_2$. We know that $a$, $c_1$, and $c_2$ make a set, so $a$ and $b_1$ cannot form a set with $c_1$ or $c_2$. The same is true for $b_2$. Thus, Bob has only one card left to block a set, and Alice has two ways to create a set using $a$ with  $b_1$ or  $a$ with $b_2$. Thus a winning card for Alice is available, and she wins on her fourth move.

\section{Our classification using groups}

Now we use the notion of groups to explain what we are doing. See, for example, \cite{C} for more a comprehensive introduction to group theory.

We start with a formal definition.

A group is a set, $G$, together with an operation $\cdot$ which is called multiplication. That is, we can multiply two elements of a group $a$ and $b$ to form another element, denoted $a \cdot  b$ or $ab$. The set and operation, $(G, \cdot)$, must satisfy four group axioms:

\begin{itemize}
\item \textbf{Closure.} For all $a$, $b$ in $G$, the result of the operation, $a\cdot b$, is also in $G$.
\item \textbf{Associativity.} For all $a$, $b$ and $c$ in $G$, $(a \cdot  b) \cdot c = a \cdot (b \cdot  c)$.
\item \textbf{Identity element.} There exists an element $e$ in $G$ such that, for every element $a$ in $G$, the equation $e \cdot a = a \cdot e = a$ holds. Such an element is unique and is called \textit{the identity element}.
\item \textbf{Inverse element.} For each $a$ in $G$, there exists an element $b$ in $G$, denoted $a^{-1}$, such that $a \cdot b = b \cdot a = e$, where $e$ is the identity element.
\end{itemize}

We are interested in the following group $G$ that acts on the cards; we want to shuffle features and values within one feature. The features can be shuffled in 24 ways. For each feature, we can shuffle the values in 6 ways. Thus, the size of our group $G$ is
\[24 \cdot 6^4 = 31104.\]

We can use group $G$ to transform any card to any other card. As mathematicians say, all cards are in the same orbit of the group action. Given a card, the subgroup of $G$ that keeps this card in place is called a \textit{stabilizer subgroup} of this card. What is the number of elements, often referred to as the order, of the stabilizer subgroup? Suppose we have one red empty oval on the card. Then, we can make any shuffle of the features. For example, if we swap color and shading as long as empty is swapped with oval, the card stays in place. In addition, if the features stay in place, for each feature we can swap two values that are not on the card. For example, we can swap green and purple. That means the stabilizer group has order $24 \cdot 2^4 = 384$.

Now we use a famous theorem from group theory that states that the number of elements in an orbit is equal to the size of the group divided by the size of the stabilizer subgroup of any of the elements in the orbit. As all cards in the deck are in the same orbit, the number of cards must be $31104/384 = 81$. We already knew that, but these calculations provide an example of how group theory might help us.

\subsection{Acting on sets}

Now we want to calculate a stabilizer of a set. We start with a set of order 1. As an example, we use the set in Figure~\ref{fig:set1}. In this set the color, the shading, and the shape are the same: the color is purple, the shading is solid, and the shape is diamonds. For each feature that has the same value for all the cards, we can swap the two values that do not occur, making $2^3$ ways in total. In our particular example, we can swap green with red, empty with striped, and ovals with squiggles. We can also shuffle these features as long as we keep the value we need. For example, we can swap color with shading as long as purple is swapped with solid. There $3! = 6$ ways to shuffle the features that are the same. Also, we can shuffle the values of the feature that is different in the set, in $3!$ ways. This is equivalent to reordering the cards themselves. So the stabilizer is $2^3\cdot 3! \cdot 3!=288$. Then by dividing we get $31104/288=108$, the number of order 1 sets.

We can generalize the example above to calculate the stabilizer of a set of order $k$. Recall that $k$ is the number of features that are different in the set. Correspondingly, $4-k$ is the number of features that are the same in the set. We can swap values that do not appear for features that are the same in the set in two ways for each feature for the total of $2^{4-k}$. We can swap features that are the same in the set with each other for the total of $(4-k)!$. We can swap features that are different in the set with each other for the total of $k!$. Finally, we can shuffle the cards in the set $3! = 6$ ways. Thus, the order of the stabilizer is 
\[6\cdot 2^{4-k} (4-k)!k!.\]

Now we explain why all the sets of the same order belong to the same orbit of the group. The cards of an order-1 set can always be written as 0000, 0001, 0002, by assigning different number values to different values of features. This is equivalent to saying that any set of order 1 belongs to the same orbit as the set written as 0000, 0001, 0002. Thus, two order 1 sets belong to the same orbit. The same can be said about order 2, 3, and 4 sets. Order 2 sets can always be written as 0000, 0011, 0022, order 3 sets can always be written as 0000, 0111, 0222, and order 4 sets can always be written as 0000, 1111, 2222.

As all the sets of the same order belong to the same orbit, the number of different sets of order $k$ is the order of the group divided by the order of the stabilizer. That is, we divide $31104 = 24 \cdot 6^4$ by $6\cdot 2^{4-k} (4-k)!k!$. The result is 
\[\binom{4}{k}3^{4-k}6^{k-1}.\]

Not surprisingly, this is the same answer we got before.

\subsection{Acting on squares}

When we use an element of this group acting on the cards, the order of each set doesn't change. That means when the group acts on the magic square the order of each triplet doesn't change. 

We can also transform any set to any other set of the same order. Can we transform any magic SET square with the given orders of triplets to any other set square with the same orders of triplets? Yes, we can. We show it using an example. Suppose the squares have type 1, 2, 3, 3. From Table~\ref{table:orderfeature}, the first square has feature distribution $x_1$, $y_1z_1$, $x_1y_1z_1$, and $x_1y_1z_1$, while the second square has feature distribution $x_2$, $y_2z_2$, $x_2y_2z_2$, and $x_2y_2z_2$. Our transformation should send features $x_1$, $y_1$, $z_1$, and $w_1$ to $x_2$, $y_2$, $z_2$, and $w_2$. After that, we can pick the values of each feature in such a way that the transformation sends the first square to the second square. We can do it for a square of any type in a similar manner.

Now we want to find the subgroups of the group that keep a particular magic SET square fixed. The group can't change the order of each triple. That means if rows and columns are of different orders we can't swap them.

Consider a magic SET square of order 2. Consider a square of type (1-1;2-2) as in Figure~\ref{fig:(1-1;2-2)}. What elements of the group keep it in place? First, we check if we can shuffle features. Columns differ only by number, so the number feature together with the values has to stay in place. Similarly, shape has to be in place. Color and shading have unique values in the square, so we can swap them as long as red is swapped with empty. In addition, while keeping features in place we can swap green and purple, and also, solid and striped. That means the stabilizer has order 8. That means the orbit is of size 3888.

In general, the procedure is as follows. Suppose $a$, $b$, $c$, and $d$ are the values of the triplets for our square of order 4. We can swap features that are the same in $(4-k)!$ ways. Also, we can swap the values that are not present in the squares for the features that are the same in $2^{4-k}$ ways. Now we look at features that are different through the square. We can swap features that are different in the same set of triplets. As each feature is different in exactly three triplets, it is equivalent to swapping features that are the same in one triplet. This contributes $(4-a)!(4-b)!(4-c)!(4-d)!$ for the total of
\[2^{4-k}(4-k)!(4-a)!(4-b)!(4-c)!(4-d)!.\]

Thus we have a formula for a square with a given arrangement of orders:
\[N(a,b,c,d) = \frac{3888\cdot 2^k}{(4-k)!(4-a)!(4-b)!(4-c)!(4-d)!}.\]

Not surprisingly, these are the same numbers we got before. Group theory helps to do it faster and neater.

\section{Acknowledgments}

This project was done as part of MIT PRIMES STEP, a program that allows students in grades 6 through 9 to try research in mathematics. Tanya Khovanova is the mentor of this project. We are grateful to PRIMES STEP and to its director, Slava Gerovitch, for this opportunity.

\end{document}